\documentclass[11pt]{amsart}

\usepackage{amsthm,amssymb,amsmath,epsfig,graphics,appendix,bm, mathrsfs,bbm}
\usepackage[all]{xy}

\usepackage{xcolor}
\usepackage[T1]{fontenc}
\usepackage{stmaryrd} 
\usepackage[left=1in, right=1in, top=1.1in,bottom=1.1in]{geometry}
\usepackage{enumitem}
\usepackage{hyperref}
\hypersetup{
     colorlinks   = true,
     citecolor    = blue,
     linkcolor=blue
}
\usepackage{cancel}
\usepackage{comment}

\usepackage{csquotes}
\usepackage{graphicx}
\usepackage{pstricks}
\usepackage{lmodern}

\newfam\msbmfam\font\tenmsbm=msbm10\textfont
\msbmfam=\tenmsbm\font\sevenmsbm=msbm7
\scriptfont\msbmfam=\sevenmsbm
\def\NN{\mathbb N}
\def\ZZ{\mathbb Z}

\def \RR{\mathbb R} \def \R{\mathbb R}

\def\E{{\bb E}}

\def\E{{\mathbb E}}

\def\proof{{\noindent \it Proof\quad  }}

\def\rme{{\rm
e}}

\def\be{{\beta}}

\def\la{{\lambda}}
\def\si{{\sigma}}

\def\Om{{\Omega}}
\def\be{{\beta}}

\def\si{{\sigma}}
\def\la{{\lambda}}

\def\i{{\imath}}

\def\bbx{{\boldsymbol{x}}}
\def\bbX{{\boldsymbol{X}}}
\def\bby{{\boldsymbol{y}}}

\def\bbz{{\boldsymbol{z}}}

\def\bbk{{\boldsymbol{k}}}
\def\bbn{{\boldsymbol{n}}}

\def\bbi{{\boldsymbol{i}}}
\def\bbt{{\boldsymbol{t}}}

\def\bbs{{\boldsymbol{s}}}\def\bbr{{\boldsymbol{r}}}
\def\bbu{{\boldsymbol{u}}}
\def\bbv{{\boldsymbol{v}}}
\def\bbl{{\boldsymbol{l}}}
\def\bbxi{{\boldsymbol{\xi}}}
\def\gg{\mathfrak g}

\def\la{\langle}
\def\ra{\rangle}
\def\I{\mathbf I}

\def\llb{\llbracket}
\def\rrb{\rrbracket}
\def\deq{\overset{\rm{def}}{=}}
\def\ep{\varepsilon}

\numberwithin{equation}{section}
\newtheorem{theorem}{Theorem}[section]
\newtheorem{thm}{Theorem}[section]
\newtheorem{lemma}[thm]{Lemma}

\newtheorem{Def}[thm]{Definition}
\newtheorem{prop}[thm]{Proposition}

\newtheorem{remark}[thm]{Remark}

\newtheorem{notation}[thm]{Notations}

\renewcommand{\theequation}{\arabic{section}.\arabic{equation}}

\def\ZZ{{\mathbb Z}}
\def\PP{{\mathbb P}}

\def\Ad2{{| A-\widetilde A |_{L^{2}_{loc}} }}

\def\om{{\omega}}\def\Om{{\Omega}}

 \renewcommand{\theequation}{\arabic{section}.\arabic{equation}}
\renewcommand{\theequation}{\arabic{section}.\arabic{equation}}
\newcommand{\D}{{\mathrm d}}

\newcommand{\nc}{\newcommand}
\nc{\bi}{\bibitem}
\newcommand{\e}{{\rm e}}
\def\wc{\overset{d}{\longrightarrow}}

\DeclareTextFontCommand{\emph}{\em}

\begin{document}
\title[Scaling limit of a long-range random walk]{Scaling limit of a long-range random walk in time-correlated random environment}

\author[G. Rang] {
Guanglin  Rang}
\address{School of Mathematics and Statistics, Wuhan University, Wuhan 430072,China;
Computational Science Hubei Key Laboratory, Wuhan University, Wuhan, 430072, China} 
\email{glrang.math@whu.edu.cn}
\author[J. Song]{Jian Song}\address{ Research Center for Mathematics and Interdisciplinary Sciences, Shandong University, Qingdao 266237, China;
School of Mathematics, Shandong University, Jinan 250100, Shandong, China}
\email{txjsong@sdu.edu.cn}

	\author[M. Wang]{Meng Wang}
	\address{School of Mathematics, Shandong University, Jinan, Shandong, 250100, China}
	\email{wangmeng22@mail.sdu.edu.cn}

\date{}
\subjclass[2010]{60F05; 60H15; 82C05}

		\keywords{Directed polymer, scaling limit, stochastic heat equation, fractional noise}

\maketitle

\tableofcontents
\begin{abstract}
This paper concerns a long-range random walk in random environment in dimension $1+1$, where the environmental disorder is independent in space but has long-range correlations in time. We prove that two types of  rescaled partition functions converge weakly to the Stratonovich solution and the It\^o-Skorohod solution respectively of a fractional stochastic heat equation with multiplicative Gaussian noise which is white in space and colored in time.
\end{abstract}



\section{Introduction}\label{sec:introduction}

The model of  directed polymer in random environment was first introduced by Huse and Henley \cite{hh85} in the study of the Ising model and  became a canonical model for disordered systems (see e.g. the lecture notes by Comets \cite{comets17}). In recent years, it has attracted much attention in particular due to its intimate connections to the stochastic heat equation (parabolic Anderson model), the stochastic  Burgers equation, and the Kardar-Parisi-Zhang (KPZ) equation and its universality class (see \cite{comets17} and \cite{corwin12} for a review). 

For a simple symmetric random walk  in i.i.d. random environment on $\NN\times \ZZ$,  Alberts et al.~\cite{alberts2014}   proved that the rescaled partition function  converges weakly to the \emph{It\^o-Skorohod} solution of the stochastic heat equation (SHE) with multiplicative space-time  Gaussian white noise. This result was extended by Caravenna et al. \cite{CAR-SUN-ZYG} to a long-range random walk in i.i.d. random environment and to other disordered models (disordered pinning model and random field Ising model). Later on, Rang \cite{rang2020} proved that, for a simple   random walk in random environment which is white in time but correlated in space,  the rescaled partition function converges weakly to the \emph{It\^o-Skorohod} solution of SHE with Gaussian noise  white in time and colored in space;  this result was extended to long-range random walks by Chen and Gao \cite{cg23}; if the random environment is given by the occupation field of a Poisson system of independent random walks on $\ZZ$ which is now correlated in both time and space,   Shen et al. \cite{Shen-Song-Sun-Xu2021} showed that the scaling limit is the \emph{Stratonovich} solution of SHE with space-time colored Gaussian noise whose covariance coincides with the heat kernel.  It is suggested by the results and methodologies in the above-mentioned papers that the \emph{temporal independence} of the random environment plays a critical  role when identifying the scaling limits for partition functions, which will be further discussed in Section~\ref{sec:main-result}.

In this paper, we aim to study the scaling limit of the partition function for a long-range random walk in  random environment on  $\NN\times \ZZ$ where the disorder is independent in space but has long-range correlations in time. 

We remark that the model of directed polymer in space-time correlated  random environment whose covariance  has a power-law decay was first considered in physics literature by Medina et al. \cite{Medina1989Burgers} where the Burgers equation with colored noise was studied and then applied to analyse directed polymer and interface growth.  We also would like to mention that Rovira and Tindel \cite{rt05} introduced Brownian polymer in a centered Gaussian field that is white in time and correlated in space on $\R_+\times\R$ and studied the asymptotic behavior of the partition function; in the two subsequent papers, Bezerra et al. \cite{btv08} obtained the superdiffusivity and  Lacoin  \cite{Lacoin2011} investigated the effect of strong spatial correlation. Finally, we remind the reader that the model of  long-range directed polymer has been studied by Comets \cite{comets07}, and more recently by Wei \cite{wei16}.

\subsection{Notations and known results}  \label{sec:review}
For the reader's convenience, we  collect the mathematical notations that will be used throughout this article. 
\begin{notation}\label{notations}
	Let $\mathbb N$ denote the set of natural numbers without 0, i.e., $\mathbb N\deq\{1, 2, \dots\}$;  
	for $N\in \NN$, $\llb N\rrb\deq\{1,2,\dots,N\}$;  for $a\in\mathbb R$, $[a]$ means the greatest integer that is not greater than $a$; $\|\cdot\|$ is used for the Euclidean norm;  let $\bbk\deq(k_1, \dots, k_d),\bbl\deq(l_1, \dots, l_d), \bbx\deq(x_1, \dots, x_d),\bby\deq (y_1,\dots, y_d)$ etc stand for vectors in $ \ZZ^d$ or $\RR^d$ depending on the context;  we use $C$ to denote a generic positive constant that may change from line to line; we say $f(x)\lesssim g(x)$ if $f(x)\le C g(x)$ for all $x$;  we write $f(x)\sim g(x)$ (as $x\to \infty$), if $\lim\limits_{x\to\infty}f(x)/g(x)=1$.  We use $\wc$ to denote the convergence in distribution (also called weak convergence) for random variables/vectors. For a random variable $X$, $\|X\|_{L^p} = \left(\E[|X|^p]\right)^{1/p}$ for $p\ge 1$. 
\end{notation}

Let $S=\{S_i, i\in\mathbb N_0\}$ be a random walk in $\ZZ$  and $\om=\{\om(i,k), (i,k)\in \NN\times\ZZ\}$ be a family of random variables independent of $S$ serving as the random environment (disorder).  We shall use $\PP_{S}$ and $\E_S$ (resp. $\PP_{\om}$ and $\E_\om$) to denote the probability and expectation  in the probability space of $S$ (resp. $\om$), respectively. The probability and expectation in the product probability space of $(S, \om)$ is denoted by $\PP$ and $\E$, respectively.

 Given $N\in\mathbb N$ and $k\in \mathbb Z$, let  $S^{(N+1,k)}=\big\{S^{(N+1,k)}_i,i \in\llb N+1\rrb\big\}$ be a \emph{backward} random walk  in $\mathbb Z$ with $S^{(N+1,k)}_{N+1}=k$, and the \emph{partition function} is defined by
	 \begin{align}\label{partition2}
		Z_\om^{(N)}(\beta,k)\deq\sum_S{\rm e}^{\be \sum_{i=1}^N\om\left(i,S^{(N+1,k)}_i\right)}\PP(S)=\E_S\left[{\rm e}^{\be \sum_{i=1}^N\om\left(i,S^{(N+1,k)}_i\right)}\right],
	\end{align}
	where $\be=1/T>0$ is the inverse temperature. We stress that the random walk is backward on the time interval $[1,N+1]$ in the sense that $S^{(N+1,k)}_{N+1}=k$ while there is no restriction on $S^{(N+1,k)}_1$.  Hence, $Z_\om^{(N)}(\beta,k)$ given by \eqref{partition2} is  indeed a  point-to-line  partition function, and it corresponds to directly a discrete version of  the solution to stochastic heat equation (see Proposition \ref{prop:FK}), which facilitates our calculations.  Throughout the rest of the article, we shall omit the superscript $(N+1,k)$ for the backward random walk to simplify the notation. 

When $S$ is a simple symmetric random walk  with i.i.d. increments and $\om(i, k)$ are i.i.d. random variables with \[\lambda(\beta)\deq \log \E[\rme^{\beta\om(i, k)}]<\infty\] for some $\beta$ sufficiently small,  Alberts et al. \cite{alberts2014} introduced the so-called \emph{intermediate disorder regime}. More precisely, if $\beta$ is scaled in the  way \[\beta\to \hat\beta_N\deq\beta N^{-\frac{1}{4}},\]  one has the following weak convergence
$$
\rme^{-N\lambda(\hat \beta_N)} Z_\om^{(N)}(\hat\beta_N, N^{1/2}x)\wc u(1,x) \text{ as } N\to \infty,
$$
where $u(t,x)$ is the It\^o-Skorohod solution of the equation
\begin{equation}\label {SHE}
	\begin{cases}
		\dfrac{\partial u(t,x)}{\partial t}=\frac12\Delta u(t,x)+\sqrt 2\beta u(t,x)\dot{W}(t,x), ~ t>0, x\in\R,\\
		u(0, x)=1.
	\end{cases}
\end{equation}
Here, $\dot W$ is space-time Gaussian white noise and  It\^o-Skorohod solution means that the product in $u(t,x)\dot W(t, x)$ is a Wick product, or equivalently the associated stochastic integral is an It\^o-Skorohod integral.

Under the same setting except that the random walk $S$ may be long-range (i.e., the increments of $S$ now lie in the domain of attraction of a $\rho$-stable law with $\rho\in(1, 2]$), the result of \cite{alberts2014} was extended in  Caravenna et al. \cite{CAR-SUN-ZYG} where  a unified framework based on the Lindeberg principle for polynomial chaos expansion was developed to study the scaling limits for pinning model, directed polymer model, and Ising model. Still using the setting of \cite{alberts2014} but assuming that the disorder $\{\om(i, k), (i, k)\in \mathbb N\times \mathbb Z\}$ is correlated in space (but still independent in time),  Rang \cite{rang2020} and Chen-Gao \cite{cg23} obtained the weak convergence of the rescaled partition function to the It\^o-Skorohod solution of \eqref{SHE} with the Gaussian noise $\dot W(t, x)$ being colored in space (still white in time).

\subsection{Model description} \label{sec:model}
Motivated by the above-mentioned works, we aim to study the weak convergence of the rescaled partition function for a long-range random walk  $\{S_n, n\in\NN\}$ in random environment on $\NN\times \ZZ$, where the disorder $\{\om(i,k), i\in \NN, k\in \ZZ\}$ is correlated in time and independent in space. 

In our model, the directed polymer (random walk) is given by $S_n=\sum_{i=0}^n Y_i$, where $\{Y_i\}_{i\in\NN_0}$ are independent and identically distributed random variables  with mean zero which have a 1-lattice distribution  belonging to the domain of attraction of a $\rho$-stable distribution  with density function denoted by $\gg_\rho(\cdot)$. We assume $1<\rho\le 2$ and
\begin{equation}\label{X-condition}
	\begin{cases}
		\PP(Y_i=k) \lesssim  |k|^{-1-\rho} \text{ for } k\in{\mathbb Z}\backslash\{0\},& \text{ if } \rho\in(1,2),\vspace{0.2cm} \\
		\E[Y_i]=0 \text{ and }  \E[Y_i^2]=1,& \text{ if } \rho=2.
	\end{cases}
\end{equation}

Let $\psi(u) \deq \E[\e^{\i u Y_i}]$ be
the characteristic function of $Y_i$. Then  the 1-lattice  distribution of $Y_i$ implies that $\psi(u)$ is periodic with period $2\pi$, and furthermore,  $|\psi(u)|<1$ for all $u\in[-\pi, \pi]\backslash\{0\}$ (see, e.g., \cite[Theorem 1.4.2]{IL71}); this property of $\psi(u)$ usually makes analysis easier. The assumption of $1$-lattice distribution on $Y_i$ also yields that the random walk $S_n=\sum_{i=1}^nY_i$ has period 1 and hence is \emph{aperiodic} . We remind the reader that a simple symmetric random walk has period 2.

Denote by $P_n(k)=\PP(S_n=k)$ for $n\in \NN, k\in\ZZ$  the probability  of $S$ being at $k$ at time $n$. Then  by the  local  limit theorem for $\rho$-stable distribution  (see \cite[Theorem 6.1]{Rvacheva1954}), the convergence
\begin{equation}\label{e:llt}
	n^{1/\rho}P_n(k)-\gg_\rho\big(k/n^{1/\rho}\big)\to 0, \text{ as } n\to\infty,
\end{equation}
holds uniformly in $k$. Noting that $\gg_\rho$ is a bounded function, \eqref{e:llt} yields 
\begin{align}\label{e:Pn-bd}
	P_n(k)\lesssim   n^{-1/\rho}, \text{ for } n\in \NN, k\in\ZZ.
\end{align}
Also under the condition of \eqref{X-condition}, we have a local deviation estimation (see, e.g., \cite[Theorem 2.6]{Berger2019}):
\begin{align}\label{e:Pn-bd'}
	P_n(k)\lesssim   n|k|^{-1-\rho}, \text{ for } n\in \NN, k\in\ZZ.
\end{align}
This together with \eqref{e:Pn-bd}, we get 
\begin{align}\label{e:Pn-bd''}
	P_n(k)\lesssim   \left(n|k|^{-1-\rho}\right) \wedge n^{-1/\rho}= n^{-1/\rho} \left(|n^{-1/\rho}k|^{-1-\rho}\wedge 1\right), \text{ for } n\in \NN, k\in\ZZ ,
\end{align}
where $a\wedge b:=\min\{a,b\}$ for $a, b\in \R$.

It is well known that 
$$
\frac{S_n}{n^{1/\rho}}\wc \xi,~~~~\mbox {as} ~~n\rightarrow\infty,
$$
where $\xi$ has  the symmetric $\rho$-stable distribution with characteristic function $\exp\{-c_\rho|\eta|^\rho\}$ for some $c_\rho>0$. Letting $\gg_\rho(t,x)$  be the density function of the corresponding $\rho$-stable process,  we have
\begin{equation*}
	\int_{\RR}\e^{\imath \eta x}\gg_\rho(t,x)\D x=\e^{-c_\rho t|\eta|^{\rho} }.
\end{equation*}
Throughout the rest of the paper we will omit the subscript $\rho$ and   use $\gg(t,x)$ exclusively  to denote the density function of the $\rho$-stable process $X$. Note that $\gg$ has the following scaling property
\begin{equation}\label{e:g-scaling}
	\gg(t,x)=t^{-\frac{1}{\rho}}\gg(1,t^{-\frac{1}{\rho}}x),
\end{equation}
and the upper bound in parallel with \eqref{e:Pn-bd''}
\begin{align}\label{e:g-bd}
	\gg(t,x)\lesssim (t|x|^{-1-\rho}) \wedge t^{-1/\rho}= t^{-1/\rho} \Big(|t^{-1/\rho}x|^{-1-\rho}\wedge 1\Big),\quad t\in \R_+,~x\in\RR.
\end{align}

We assume that the disorder in the environment is given by a family of Gaussian random variables  $\{\om(n,k),n\in\NN, k\in \ZZ\}$ with mean zero and covariance
\begin{equation}\label{e:con-gamma}
	\E[\om(n,k)\om(n',k')] = \gamma(n-n')\delta_{kk'},
\end{equation}
where $\delta_{kk'}$ is the Kronecker delta
function, i.e., $\delta_{kk'}=1$ if $k=k'$ and $\delta_{kk'}=0$ otherwise, and $\gamma(n)$ has a power law decay:
\begin{equation}\label{e:gamma-bd}
	\gamma(n)\lesssim  |n|^{2H-2}\wedge 1, \text{ for } n\in \ZZ,
\end{equation}
where  $H \in(1/2, 1]$ denotes a fixed constant throughout this paper.  We further assume that  for all $t\in\R\backslash\{0\}$,
\begin{equation}\label{e:gamma-lim}
	\lim_{N\to\infty}N^{2-2H}\gamma([Nt]) =|t|^{2H-2}. 
\end{equation}

\subsection{Main result, strategy and discussions} \label{sec:main-result}

Consider the stochastic fractional heat equation on $\R$
\begin{equation}\label {SHE1}
	\begin{cases}
		\dfrac{\partial u(t,x)}{\partial t}= -c_\rho(-\Delta)^{\frac{\rho}{2}} u(t,x)+\beta u(t,x)\dot{W}(t,x), \\ 
		u(0,x)=1,
	\end{cases}
\end{equation}
where $\rho\in(1,2]$ and $\dot{W}(t,x)$ is Gaussian noise  with covariance function given by
\begin{equation}\label{e:cov}
	\E\bigg[\dot {W}(t,x)\dot {W}(s,y)\bigg]= K(t-s, x-y)\deq |t-s|^{2H-2}\delta(x-y), 
\end{equation}
with $H\in(1/2, 1]$ and $\delta(\cdot)$ being the Dirac delta function. In particular, $\dot W$ is a spatial white noise independent of time if $H=1$.  

We consider two types of solutions of \eqref{SHE1}---the \emph{Stratonovich solution} if the product  $u(t, x) \dot W(t, x)$ is an ordinary product (i.e., the associated stochastic integral is a Stratonovich integral) and the \emph{It\^o-Skorohod solution} if $u(t, x) \dot W(t, x)$ is  a Wick product (i.e., the associated stochastic integral is a Skorohod integral). See Section \ref{sec:malliavin} for details.  

Now we are ready to present our main results.  For the first type of partition function $Z_\om^{(N)}$ which is given in \eqref{partition2}, we have the following result. 
\begin{theorem}\label{thm:main}
	Let the random walk $S$ and the disorder $\om$ be given as in Section \ref{sec:model}.    Let $H$ and $\rho$ be parameters satisfying 
	\begin{equation}\label{e:con}
		H\in(1/2, 1]~ \text{ and }  ~\theta\deq H-\frac1{2\rho}>\frac12.
	\end{equation}
	Consider  the rescaled  partition function $Z_\om^{(N)}(\hat\beta_N,  N^{1/\rho}x_0)$ given in \eqref{partition2} under the scaling \[\beta\to \hat\beta_N\deq\beta N^{-\theta}.\]
	Let $u(t, x)$ be the \emph{Stratonovich solution} of \eqref{SHE1}.   Then we have
	\[Z_\om^{(N)}(\hat\beta_N,  N^{1/\rho} x_0) \wc u(1, x_0), \text{ as } N\to \infty.\]
\end{theorem}

In comparison with \eqref{partition2}, another type of point-to-line partition function is  given by
\begin{align}\label{partition3}
	\tilde Z_\om^{(N)}(\beta,k)\deq \E_S\left[{\rm e}^{\be  \sum_{i=1}^N\om(i,S_i) -\frac{\beta^2}{2}\sum_{i,j=1}^N \gamma(i-j) \mathbf 1_{\{S_i=S_j\}}}\right],
\end{align}
where   $S=S^{(N+1,k)}=\big\{S_i,i \in\llbracket N+1\rrbracket\big\}$ is a backward random walk with $S_{N+1}=k$. The extra term in the exponential is half of the variance of $\be  \sum_{i=1}^N\om(i,S_i)$ conditional on the random walk $S$. Indeed, we have
\begin{equation}\label{e:correction}
	\begin{aligned}
		\E_\om\left[\left(\sum_{i=1}^N \om(i,S_i) \right)^2 \right]=&
		\E_\om\left[\sum_{i=1}^N \sum_{k\in\ZZ}\om(i,k) \mathbf 1_{\{S_i=k\}}\sum_{j=1}^N \sum_{l\in\ZZ}\om(j, l) \mathbf 1_{\{S_j=l\}}\right]\\
		=&\sum_{i,j=1}^N \sum_{k,l\in\ZZ}\gamma(i-j) \delta_{kl} \mathbf 1_{\{S_i=k\}}\mathbf 1_{\{S_j=l\}}\\
		=&\sum_{i,j=1}^N \gamma(i-j) \mathbf 1_{\{S_i=S_j\}},
	\end{aligned}
\end{equation}
which can be viewed as a weighted intersection local time of $S$.

For the scaling limit of $\tilde Z_\om^{(N)}$ given in \eqref{partition3}, we have the following result in parallel with Theorem~\ref{thm:main}.
\begin{theorem}\label{thm:main1}
	Let the random walk $S$ and the disorder $\om$ be given as in Section \ref{sec:model}.    Let $H$ and $\rho$ be parameters satisfying 
	\begin{equation}\label{e:con'}
		H\in(1/2, 1],~ \rho \in (1,2], \text{ and }  ~\theta= H-\frac1{2\rho}>0.
	\end{equation}
	Consider  the rescaled  partition function $\tilde Z_\om^{(N)}(\hat\beta_N,  N^{1/\rho}x_0)$ given in \eqref{partition3} under the scaling \[\beta\to \hat\beta_N\deq\beta N^{-\theta}.\]
	Let $\tilde u(t, x)$ be the \emph{It\^o-Skorohod solution} of \eqref{SHE1}.   Then we have
	\[\tilde Z_\om^{(N)}(\hat\beta_N,  N^{1/\rho} x_0) \wc \tilde u(1, x_0), \text{ as } N\to \infty.\]
\end{theorem}
\begin{remark}\label{rem:rho}
	In Theorem \ref{thm:main1}, the condition $\rho>1$ in \eqref{e:con'} arises to guarantee the existence and uniqueness of the It\^o-Skorohod solution of \eqref{SHE1} (see Remark \ref{rem:sko-sol} in Section \ref{sec:she1}). In Theorem~\ref{thm:main}, the condition \eqref{e:con} also implies $\rho>\frac1{2H-1}\ge 1$.  This is why we restrict us on the case $\rho\in(1, 2]$ in this article. 
\end{remark}
\begin{remark}\label{rem:con-two-cases}
	For \eqref{SHE1},  it requires less restrictive conditions to have an It\^o-Skorohod solution than to have a Stratonovich solution (see Remark \ref{rem:con-two-solutions} and see also \cite{Song2017} for a more general class of SPDEs).  This explains why condition \eqref{e:con'} for the It\^o-Skorohod case is  weaker than condition \eqref{e:con}  for the Stratonovich case.
\end{remark}

The proofs of Theorems \ref{thm:main}  and \ref{thm:main1}   will be presented in Section \ref{sec:wc-partition}. For the reader's convenience, here we explain briefly the strategy  for Theorem \ref{thm:main} (the proof of Theorem~\ref{thm:main1} is similar but easier).  By Taylor's expansion, the rescaled partition function can be written as 
\begin{align*}
	Z_\om^{(N)}(\hat\beta_N,N^{1/\rho} x_0)&=\sum_{m=0}^\infty\frac{1}{m!}{\mathbb S}_m^{(N)},
\end{align*}
where
\begin{align}\label{e:bbS-m}
	{{\mathbb S}}_m^{(N)}=\hat\beta^m_N\sum_{n_1, \dots, n_m\in \llb N \rrb}\sum_{ k_1,\dots, k_m\in \ZZ}\om(n_1, k_1)\cdots\om(n_m, k_m)P_{\bbn^*}(N^{1/\rho} x_0; k_1,\dots, k_m),
\end{align} 
with $P_{\bbn^*}$ (see eq. \eqref{e:P*}) being the product of the transition densities of the random walk $S$. Meanwhile, the Stratonovich solution of the continuum equation \eqref{SHE1} has the following series representation (see eq. \eqref{chaos-solution}):
\[u(1, x_0) =\sum_{m=0}^\infty \beta^m \mathbb I_m (\gg_m (\cdot;1, x_0)),\]
where $\mathbb I_m(\cdot)$ is a  multiple Stratonovich integral (see Section \ref{sec:Malliavin}) and $\gg_m$ (see eq. \eqref{e:gm}) is the product of the transition densities of the stable process $X$.  We remind the reader that the continuum  multiple Stratonovich integral   $\mathbb I_m(\gg_m(\cdot;1,x_0))$ resembles the discrete sum $\mathbb S_m^{(N)}$.

Under condition \eqref{e:con}, one can obtain the $L^1$-convergence of the series of $u(t,x)$(see Proposition~\ref{prop:L2-conv-spde}) and the uniform (in $N$) $L^1$-convergence of the series of $Z_\om^{(N)}(\hat \beta_N, N^{1/\rho} x_0)$ (see \eqref{e:uniform-convergence}). Thus, in order to prove the weak convergence $Z_\om^{(N)}(\hat\beta_N,  N^{1/\rho} x_0) \wc u(1,  x_0)$, in light of Lemma~\ref{lem:weak-con}, it suffices to obtain joint weak convergences for $\mathbb S_m^{(N)}$'s, that is,  for $k\in\NN$ and $l_1, \dots, l_k\in \NN$,
\begin{equation*}
	\left(\frac{1}{l_1!}\mathbb S_{l_1}^{(N)}, \dots, \frac{1}{l_k!}\mathbb S_{l_k}^{(N)}\right)  \wc \Big(\beta^{l_1}\mathbb I_{l_1}(\gg_{l_1}(\cdot;1,x_0)), \dots, \beta^{l_k}\mathbb I_{l_k}(\gg_{l_k}(\cdot;1,x_0))\Big), \text{ as } N\to \infty,
\end{equation*}
which is proved in Proposition \ref{m-thm}.  

Multiple Wiener integrals $\I_m(f)$ and multiple Stratonovich integrals $\mathbb I_m(f)$ are linked via the celebrated Hu-Meyer's formula \eqref{e:hm}. Usually it is more convenient to deal with multiple Wiener integrals $\I_m(f)$ whose second moment is easier to calculate.   For this reason,  in Section \ref{sec:U-stat} we first consider $U$-statistics given in \eqref{U-stat-I} which are multi-linear Wick polynomials of $\om(n_i,  k_i)$'s, and prove that they converge weakly to multiple Wiener integrals (see Proposition \ref{prop:m-integral-I}). This together with  Hu-Meyer's formula \eqref{e:hm} then  yields Proposition \ref{m-thm}.   

Note that the exponential term in the  partition function $\tilde Z_\om^{(N)}$ in \eqref{partition3} is actually a \emph{Wick exponential} (see eq. \eqref{e:wick-exp}) of $\beta\sum_{i=1}^N \om(i, S_i)$ conditional on the random walk $S$. Thus, to prove Theorem \ref{thm:main1}, one only needs Proposition \ref{prop:m-integral-I} on the weak convergence of multi-linear Wick polynomials.

As can be seen, our approach is very much inspired by \cite{alberts2014} where the  random variables of disorder are i.i.d.  However, in our setting the disorder  is correlated in time, and it turns out that this difference is critical in the sense that it prevents us from using the trick of ``modified partition function'' (see  in \cite[eq. (4)]{alberts2014},  which was also employed in \cite{rang2020}), as this method relies on  the temporal independence of the disorder. For the same reason, the Mayer expansion (see \cite[equations (1.5) and (1.7)]{CAR-SUN-ZYG}), which was employed in \cite{CAR-SUN-ZYG} to linearize the exponential while introducing an error term that comes from the It\^o correction, cannot be applied here either, since it works well only if the disorder is independent in  both time and space.

To conclude the introduction, we make some remarks on our main results Theorems \ref{thm:main} and  \ref{thm:main1}. 

\noindent{\bf (i)} In the polymer model, the  temporal independence of the disorder plays a critical  role. Recently,  a directed polymer in time-space correlated random environment was studied in  \cite{Shen-Song-Sun-Xu2021} and it was shown that the rescaled partition function converges weakly to a Stratonovich solution.  Assuming that the disorder is correlated in time and white in space in our setting, the  rescaled partition functions  $Z_\om^{(N)}$ and $\tilde Z_\om^{(N)}$  converge weakly to  \emph{Stratonovich} and \emph{It\^o-Skorohod} solutions, respectively.  In contrast, if the disorder possesses temporal independence, after a proper scaling, the limit of the rescaled partition function is an \emph{It\^o-Skorohod} solution but not a \emph{Stratonovich} solution, as has been shown in \cite{alberts2014, CAR-SUN-ZYG, rang2020, cg23}. Among others, a technical explanation for the critical role of the temporal independence of the disorder  is the following: when the temporal independence of disorder appears, by using the ``modified partition function'' trick, the term corresponding to $\mathbb S_m^{(N)}$ given by \eqref{e:bbS-m} becomes a multi-linear polynomial of $\om(n_i,  k_i)$ for $i=1,\dots, m$ with $n_1<n_2<\cdots < n_m$, which  converges weakly to a \emph{multiple Wiener integral} due to the independence (see \cite{moo10, CAR-SUN-ZYG}).

\noindent{\bf (ii)} In contrast to the case that the disorder is independent in time where one can employ the ``modified partition function'' trick (or Mayer's expansion for the case that the disorder is independent both in time and space) to make  time indices in the summation distinct from each other, we expand the partition function directly via Taylor's expansion, and as a consequence the time indices $n_1, \dots, n_m$ in the summation \eqref{e:bbS-m} can be repeated. It turns out that  this in fact gives negligible contribution for the Stratonovich case, due to the condition $\theta>\frac12$ in  \eqref{e:con} which is more restrictive than the condition $\theta>0$ in  \eqref{e:con'} for the Skorohod case (see also \cite{alberts2014, CAR-SUN-ZYG, rang2020}). For instance, when $m=2$, the expectation of the  sum of the diagonal terms in \eqref{e:bbS-m} is $\beta\gamma(0) N^{1-2\theta}$ which converges to zero as $N$ tends to infinity if we assume  $\theta >1/2$. 

\noindent{\bf (iii)}  The study of the scaling limit of partition functions was initiated in \cite{alberts2014} in order to understand the polymer behavior in the so-called \emph{intermediate disorder regime} which sits between weak and strong disorder regimes. Meanwhile, as pointed out in \cite{CAR-SUN-ZYG}, the fact that the  rescaled partition functions converge
 weakly to a non-trivial limit indicates that the directed polymer model is \emph{disorder relevant}, since it implies that the presence of  disorder, no matter how small it is,  changes the qualitative features of the underlying homogeneous model.

\noindent{\bf (iv)}  As pointed out in \cite{Medina1989Burgers}, it is of interest to study directed polymer in random environment which has long-range correlations in time and space. In light of the work of \cite{Shen-Song-Sun-Xu2021}, we expect that our results and methodology can be extended to the model in dimension $1+d$ with space-time long-range  correlated disorder. Note that the corresponding continuum SPDEs have been investigated in \cite{HNS2011, Song2017}.
For instance, consider the limiting Gaussian noise with covariance  \[\E\bigg[\dot {W}(t,x)\dot {W}(s,y)\bigg]=|t-s|^{-\alpha_0}|x-y|^{-\alpha},\]  
where $\alpha_0\in(0,1)$ and $\alpha\in(0,d)$. Then \eqref{SHE} has a Stratonovich solution if  $\alpha<\rho(1-\alpha_0)$ (see Remark~3.1 in \cite{Song2017}) and an It\^o-Skorohod solution if $\alpha < \rho$ (see Remark~5.1 and Theorem 5.3 in \cite{Song2017}). In particular, this allows to consider the case $d>1, \rho\in(0, 1]$ (recall that the spatial independence forces us to consider only the case $d=1, \rho\in(1,2]$, see Remarks \ref{rem:rho} and \ref{rem:sko-sol}).

\noindent{\bf (v)} In our model, we assume the Gaussianity  of the disorder $\om$ for technical reasons. The loss of Gaussian property in $\om$ shall cause a lot  more complexity in computations (see e.g. \eqref{e:E-wick-prod} and \eqref{e:E-wick-prod'} for a comparison of the computations of moments for Gaussian and non-Gaussian random variables). The functional analytic approach developed in \cite{Shen-Song-Sun-Xu2021} might be a  choice to circumvent this difficulty.   Nevertheless, we conjecture that our result still holds for sub-Gaussian disorder.

\noindent{\bf(vi)} Assume the disorder $\{\om(i,k), i\in\NN, k\in\ZZ\}$ is a family of i.i.d. standard Gaussian random variables.  In this situation,  the exponent in the correction term of  the rescaled partition function $\tilde Z_\om^{(N)}(\hat \beta_N, k)$ becomes, recalling  \eqref{e:correction} and replacing $\gamma(i-j)$ by $\delta_{ij}$ therein,  \[\frac12{\hat\beta_N^2}\sum_{i,j=1}^N \delta_{ij}\mathbf 1_{\{S_i=S_j\}}=\frac12\hat \beta_N^2 N=\frac{\beta^2}2 N^{\frac1\rho},\]
which is now independent of the random walk $S$.  Thus, the rescaled  partition function for the It\^o-Skorohod case can be written as 
\[\tilde Z_\om^{(N)}(\hat \beta_N, k) = {\rm e}^{-\frac{\beta^2}2N^{\frac1\rho}}Z_\om^{(N)}(\hat \beta_N, k),\]
which converges weakly to the It\^o-Skorohod solution of \eqref{SHE1} with space-time white noise. 
This is consistent with  \cite[Theorem 2.1]{alberts2014},  noting that $N\log \E[{\rm e}^{\hat \beta_N \om(i,k)}]= \frac12\hat \beta_N^2 N$ if we assume the disorder is Gaussian.

\noindent{\bf(vii)} If the disorder $\om$ and the limiting noise $\dot W$ live in the same probability space  such that  $\dot W$ is the scaling limit of $\om$,  we can get  strong convergence (a.s. convergence) in Theorem~\ref{thm:main} and Theorem \ref{thm:main1}.  In this sense,  the rescaled partition functions $Z_\om^{(N)}$ and $\tilde Z_\om^{(N)}$ can be viewed as approximations for the Stratonovich solution and the It\^o-Skorohod solution of equation~\eqref{SHE1}, respectively. See e.g. Foondun at el.  \cite{fjl18}  and  Joseph et al.  \cite{jkm17} for related results. 

\noindent{\bf(viii)} As a byproduct, we obtain a Jensen type inequality for the integrals induced by fractional Brownian motion (see Lemmas \ref{lem:Jensen} and \ref{lem:jensen}) which is new in the literature to our best knowledge.  We provide a new approach to prove  the exponential integrability of the weighted self-intersection local time of the $\rho$-stable process in Proposition \ref{prop:exp-int}  which plays a critical role in obtaining the Feynman-Kac formula for the solution of \eqref{SHE1}; this approach also simplifies the proof  in \cite{HNS2011}  which works exclusively for the weighted self-intersection local time of Brownian motion.

The rest of the paper is organized as follows. In Section \ref{sec:SPDE}, we provide some preliminaries on Gaussian spaces and then study the limiting SPDE \eqref{SHE1}. The main results Theorem~\ref{thm:main} and Theorem \ref{thm:main1}  are proved in Section~\ref{sec:proof}. In Appendices  \ref{sec:wick}, \ref{sec:appB}, and \ref{sec:appC}, we collect some preliminaries on Wick products,  convergence of probability measures, and some other miscellaneous results that are used in this article. 

 \section{On the SPDE}\label{sec:SPDE}
In this section, we first recall some preliminaries on Gaussian spaces and then we study the limiting continuum SPDE \eqref{SHE1}.

\subsection{Preliminaries on Gaussian spaces}\label{sec:malliavin}

In this subsection, we  provide some preliminaries on Gaussian spaces. We refer to \cite{ Hu16, janson1997, nualart06} for more details.

\subsubsection{Banach spaces associated with $\dot W$} On probability space $(\Om,{\mathcal F}, P)$ satisfying usual conditions, let $\dot W=\{\dot W(t,x):t\in[0,1],x\in\mathbb R\}$ be  real-valued centered Gaussian noise with covariance
\begin{align}\label{Gaussian-noise}
	\E[\dot W(t,x)\dot W(s,y)]=&K(t-s, x-y)=|t-s|^{2H-2}\delta(x-y),
\end{align}
where $\delta$ is the Dirac delta function.  The Hilbert space $\mathscr H$ associated with $\dot W$  is the completion of smooth functions with compact support under the inner product
\begin{equation}\label{e:inner}
	\begin{aligned}
		\langle f, g \rangle_{\mathscr H}=&\int_0^1\int_0^1\int_{\RR}\int_{\RR} f(s,x)K(s-t,x-y)g(t,y)\D  x\D y\D s\D t\\
		=&\int_0^1 \int_0^1 \int_{\R} |s-t|^{2H-2} f(s,x) g(t, x) \D x  \D s \D t.
	\end{aligned}
\end{equation}
We assume $H\in(1/2, 1].$  For the case $H=1$, the noise $\dot W$ is indeed  a spatial white noise which does not depend on time.

We remark that the Hilbert space ${\mathscr H}$  contains distributions (see \cite{pt00}). For our purpose, it suffices to just consider classical measurable functions. We introduce the following Banach space $\mathcal B$ which is a subset of the set $\mathscr B([0,1]\times \RR)$ of measurable functions on $[0,1]\times \R$. 

\begin{Def} We define
	\begin{equation}\label{eq-B}
		\begin{aligned}
			\mathcal B\deq \Bigg\{&f\in\mathscr B([0,1]\times \RR): \\
			&\qquad \|f\|_{\mathcal B}\deq\left(\int_0^1 \int_0^1 \int_{\R} |s-t|^{2H-2} |f(s,x)||f(t, x)| \D x  \D s \D t\right)^{\frac12}<\infty\Bigg\}.
		\end{aligned}
	\end{equation}
\end{Def}
Clearly $\mathcal B$ is a dense subset of $\mathscr H$, and if $f\in \mathscr H$ is a nonnegative function, $f$ also belongs to $\mathcal B$.

For $N\in\NN$ and $f\in\mathcal B^{\otimes m}$, we denote by  $\mathcal A_N(f)$ the conditional expectation of $f$ with respect to the $\sigma$-algebra $\mathscr B_N^m=\sigma({\mathfrak R}_N^m)$ generated by the set ${\mathfrak R}_N^m$ of rectangles 
 \begin{equation}\label{e:rectangle}
		{\mathfrak R}_N^m:=\left\{\left[\frac{{\bbi}}{N},\frac{{\bbi}+\overrightarrow{1}_m}{N}\right)\times\left[\frac{{ \bbk}}{N^{1/\rho}},\frac{{\bbk}+{\overrightarrow 1}_{m}}{N^{1/\rho}}\right):{\bbi}\in \llb N\rrb^m, { \bbk}\in {\ZZ}^{m}\right\},
\end{equation}
where ${\overrightarrow 1}_m$ is the $m$-dimensional vector of  all ones. That is, $\mathcal A_N(f)$  is defined by the average values of $f$ on the blocks $B\in {\mathfrak R}_N^m$:
\begin{align}\label{con-exp}
	{\mathcal A_N(f)}({\bbt},{\bbx})= \frac{1}{|B|}\int_{B}f({\bbs},{\bby})\D {\bbs}\D{\bby}\cdot \mathbf 1_{B}(\bbt, \bbx),
\end{align}
where $\mathbf 1_{B}$ is an indicator function and $|B|$ is the Lebesgue measure of 
$B$. 

The following Jensen type of inequality will be used to prove the weak convergence of U-statistics in Section \ref{sec:U-stat}. 
\begin{lemma}\label{lem:Jensen}
	For  $m, N\in \mathbb N$, consider $f\in\mathcal B^{\otimes m}$ and let ${\mathcal A_N(f)}$ be given by \eqref{con-exp}. Then,  we have 
	\begin{equation}\label{e:Jensen}
		\|{\mathcal A_N(f)}\|_{\mathcal B^{\otimes m}}\leq C^m \|f\|_{\mathcal B^{\otimes m}},
	\end{equation}
	where $C$ is a constant  depending only on $H$. \end{lemma}
\proof We prove \eqref{e:Jensen} for $m=1$ and the general case can be proved in a similar way. Note that by Fubini's theorem, $\mathcal A_N(f)=\mathcal A^{\text{t}}_N(\mathcal A^{\text{x}}_N(f))=\mathcal A^{\text{x}}_N(\mathcal A^{\text{t}}_N(f))$, where $\mathcal A^{\text{t}}_N$ (resp. $\mathcal A^{\text{x}}_N$) means taking average in time (resp. space) only.  Therefore, we have
\[ \|\mathcal A_N(f) \|_{\mathcal B}=\|\mathcal A_N^{\text{t}}(\mathcal A^{\text{x}}_N (f)) \|_{\mathcal B} \le C \|\mathcal A^{\text{x}}_N(f) \|_{\mathcal B},\]
for some $C$ depending on $H$ only,  where the inequality follows from  Lemma \ref{lem:jensen}.  

Now, it suffices to prove $\|\mathcal A^{\text{x}}_N(f) \|_{\mathcal B} \le \|f \|_{\mathcal B}$.  Using the identity, for $H\in(\frac12, 1)$,
\[|s-t|^{2H-2}= c_H \int_\R |s-\tau|^{H-\frac32} |t-\tau|^{H-\frac32} \D \tau\]
where $c_H$ is a finite number only depending on $H$,
we have
\begin{align*}
	\|\mathcal A^{\text{x}}_N(f) \|_{\mathcal B}^2&=  c_H \int_\R\D x \int_0^1 \D s \int_0^1 \D t \int_{\R} \D \tau |s-\tau|^{H-\frac32}|t-\tau|^{H-\frac32} \left|\mathcal A^{\text{x}}_N(f)\right|(s,x) \left|\mathcal A^{\text{x}}_N(f)\right|(t,x)\\
	&= c_H \int_\R\D x  \int_\R\D \tau \left( \int_0^1 \D s |s-\tau|^{H-\frac32}\left|\mathcal A^{\text{x}}_N(f)\right|(s,x)  \right)^2 \\
	&\le c_H \int_\R\D x  \int_\R\D \tau \left( \mathcal A^{\text{x}}_N\left(\int_0^1 \D s |s-\tau|^{H-\frac32}|f|(s,\cdot)\right)(\tau, x) \right)^2 \\
	&\le c_H \int_\R\D x  \int_\R\D \tau \left(\int_0^1 \D s |s-\tau|^{H-\frac32}|f|(s,x)\right)^2 \\
	&=\|f\|_{\mathcal B}^2,
\end{align*}
where the last inequality follows from the classical Jensen's inequality. 
\qed

\subsubsection{Chaos expansion, Wick product, multiple integrals, etc} \label{sec:Malliavin}

Recall that $\mathscr H$ is the Hilbert space associated with Gaussian noise  $\dot W$ with covariance \eqref{Gaussian-noise}. Let $\{W(f), f\in \mathscr H\}$ be an  isonormal Gaussian process with covariance 
\begin{equation}
	\E[W(f) W(g)] =\langle f, g \rangle_{\mathscr H}.
\end{equation}
In particular, if $f=\mathbf 1_{[0,t]\times[0,x]}$, we denote
$$
W(t,x)\deq W(\mathbf 1_{[0,t]\times[0,x]}).
$$
Then the fractional noise  $\dot{W}(t,x)$ can be viewed as the partial derivative $\frac{\partial^{2}}{\partial t \partial x} W(t,x)$ in the sense of distribution. 
For $f\in\mathscr H$, we also use the integral form to denote the Wiener integral $W(f)$:
\[
\int_0^1\int_{\RR}f(s,y)W(\D s,\D y)\deq W(f). 
\]

For $m\in \NN\cup\{0\}$, let
\begin{equation}\label{e:hermit}
	H_m(x)\deq (-1)^m\e^{x^2/2}\frac{\D ^m}{\D x^m}\e^{-x^2/2},\quad x\in{\mathbb R}
\end{equation}
be the $m$th {\it Hermite polynomial}. For  $g\in \mathscr H$,  the \emph{multiple Wiener integral} of $g^{\otimes m}\in \mathscr H^{\otimes m}$ can be defined via  (see e.g. \cite{hy09,nualart06})
\begin{equation}\label{e:I-m}
	\I_m(g^{\otimes m})\deq \|g\|_{\mathscr H}^m\, H_m\left(W(g) \|g\|_{\mathscr H}^{-1}\right).
\end{equation}
In particular, we have $W(g)= \I_1(g)$ and $\I_m(g^{\otimes m})=:\!W(g)^{m}\!:$, where $:\!W(g)^{m}\!:$  is a physical Wick product (see Section \ref{sec:wick}).

For $f\in \mathscr H^{\otimes m}$, let $\hat f$ be its symmetrization, i.e., 
\[\hat f(t_1, x_1, \dots, t_m,x_m) = \frac1{m!}\sum_{\sigma\in \mathcal P_m} f(t_{\sigma(1)}, x_{\sigma(1)}, \dots,t_{\sigma(m)}, x_{\sigma(m)}),\]
where $\mathcal P_m$ is the set of all permutations of $\llb m\rrb=\{1,2,\dots,m\}.$  Let ${\mathscr H}^{\hat\otimes m}$ be the symmetrization  of ${\mathscr H}^{\otimes m}$. Then for $f\in {\mathscr H}^{\hat\otimes m}$,  one can define the $m$th multiple Wiener integral $\I_m(f)$ via \eqref{e:I-m} by a limiting argument.  Moreover, for $f \in {\mathscr H}^{\hat\otimes m}$ and  $g\in {\mathscr H}^{\hat\otimes n}$, we have
\begin{equation}\label{e:ImIn}
	{\mathbb E}[\I_m(f)\I_n(g)]=m!\langle f,g\rangle_{\mathscr H^{\otimes m}}\delta_{mn}, 
\end{equation}
where we recall that $\delta_{mn}$ is the Kronecker delta function.  For $f\in {\mathscr H}^{\otimes m}$ which is not necessarily symmetric, we simply let $\I_m(f)\deq \I_m(\hat f).$   We also take the following notation for multiple Wiener integrals:
\begin{align}\label{e:multiple-Wiener}
	\int_{([0,1]\times {\mathbb R})^m}f({\bbt},{\bbx}) W(\D {t_1}, \D {x_1})\diamond\cdots \diamond W(\D {t_m}, \D {x_m})\deq \I_m(f), \text{ for } f\in\mathscr H^{\otimes m}.
\end{align}

For $f\in \mathscr H^{\hat\otimes m}$ and  $g\in\mathscr H^{\hat \otimes n}$,  their contraction of $r$ indices for $1\leq r\leq m\land n$ is defined by
\begin{align*}
	&(f\otimes_rg)(t_1,x_1,\dots, t_{m+n-2r},x_{m+n-2r})\nonumber\\
	&\deq \int_{[0,1]^{2r}}\int_{{\mathbb R}^{2r}}f(t_1,x_1,\dots,t_{n-r},x_{n-r}, u_1, y_1, \dots, u_r,y_r) \prod_{i=1}^rK(u_i-v_i, y_i-z_i)\\
	&\qquad\qquad\qquad  \times g(t_{n-r+1},x_{n-r+1},\dots, t_{n+m-2r},x_{n+m-2r}, v_1, z_1, \dots, v_r, z_r)\D{\bby} \D{\bbz}\D \bbu \D \bbv,
\end{align*}
and we have the following recursive formulas:
\begin{align}\label{contract1}
	\I_m(f)\I_n(g)=\sum_{r=0}^{m\land n}r!\binom{n}{r}\binom{m}{r}\I_{m+n-2r}(f\otimes_rg).
\end{align}
In particular, when $n=1$ we have,
\begin{align}\label{contract2}
	\I_m(f)\I_1(g)=\I_{m+1}(f\otimes g)+m\I_{m-1}(f\otimes_1g).
\end{align}

Square integrable random variables have a unique chaos expansion as stated below. 
\begin{prop}
	Let $\mathcal G$ be the $\sigma$-field generated by $W$. Then for any $F\in L^2(\Om,\mathcal G,P)$, it admits a unique chaos expansion
	$$
	F=\sum_{m=0}^\infty \I_m(f_m) ~ \text{ with } ~f_m\in \mathscr H^{\hat\otimes m},
	$$
	where the series converges in $L^2$. Moreover,
	$$
	\E [F^2]=\sum_{m=0}^\infty {m!}\|f_m\|^2_{\mathscr H^{\otimes m}}.
	$$
\end{prop}

For  $f\in \mathscr H^{\otimes m}$ and $g\in \mathscr H^{\otimes n}$, the {\it probabilistic Wick product} (or \emph{Wick product}) of $\I_m(f)$ and $\I_n(g)$ is defined by
\begin{equation}\label{e:Wick-product}
	\I_m(f) \diamond \I_n(g) \deq \I_{m+n}(f\otimes g).
\end{equation}
Unlike the \emph{physical} Wick product, the Wick product defined by \eqref{e:Wick-product}  has the same properties as the ordinary product.  As an example, by  \eqref{e:Wick-product} and \eqref{contract2}, we can show that for $f\in \mathscr H^{\otimes m}$ and $g\in \mathscr H$,
\begin{equation*}\label{e:Wick1}
	\I_m(f) \diamond \I_1(g) =\I_m(f) \I_1(g) - m\I_{m-1} (f\otimes_1g).
\end{equation*}
For two square integrable random variables  $F=\sum_{m=0}^\infty \I_m(f_m)$ and $G=\sum_{n=0}^\infty \I_n(g_n)$, their Wick product is given by 
\[F\diamond G= \sum_{m=0}^\infty\sum_{n=0}^\infty \I_{m+n}(f_m \otimes g_n),\]
as long as the right-hand side is well defined. For instance, 
let $\varepsilon(u)$ be the \emph{exponential vector} of $u$ for $u\in \mathscr H$: \[\varepsilon(u)\deq \exp\left(W(u)-\frac{\|u\|^2_{\mathscr H}}{2}\right).\] 
Then  $\varepsilon (u)\diamond\varepsilon(v)=\varepsilon(u+v)$ for $u,v\in\mathscr H$. We also recall that for a centered Gaussian random variable $F$, its \emph{Wick exponential} is given by 
\begin{equation}\label{e:wick-exp}
	\exp^\diamond(F)\deq \exp\left(F-\frac12 \E[F^2]\right)=\sum_{m=0}^\infty \frac1{m!} F^{\diamond m}. 
\end{equation}

As mentioned in Section \ref{sec:wick}, the physical and probabilistic Wick products of a Gaussian vector coincide. In particular, for $f_i\in \mathscr H, i=1, \dots, n$, we  have 
\begin{equation}\label{e:link-wick}
	:\!\I_1(f_1)\cdots \I_1(f_n)\!: = \I_1(f_1)\diamond \cdots \diamond \I_1(f_n). 
\end{equation}

In contrast to the \emph{multiple Wiener integral} given in \eqref{e:I-m}, the \emph{multiple Stratonovich integral} of $g^{\otimes m}$ for $g\in \mathscr H$ is defined by (\cite{hm, Hu16}): 
\begin{equation}\label{e:II-m}
	\mathbb I_m(g^{\otimes m}) \deq  W(g)^m. 
\end{equation}
By a limiting argument, one can define $\mathbb I_m(g)$ for $g$ in some appropriate space. We also take the following notation for $\mathbb I_m(f)$ (in comparison with $\I_m(f)$ in \eqref{e:multiple-Wiener}):
\begin{align*}
	\mathbb I_m(f)=&\int_{([0,1]\times {\mathbb R})^m}f({\bbt},\bbx) W(\D {t_1}, \D {x_1})\cdots W(\D {t_m}, \D {x_m})\\
	=&\int_{([0,1]\times {\mathbb R})^m}f({\bbt},{\bbx})\prod_{i=1}^m W(\D {t_i}, \D {x_i}). 
\end{align*}

For $f\in\mathscr H^{\otimes m}$, define the $k$th order trace $\mathrm{Tr}^kf$ of $f$ by
\begin{align}\label{e:trace-k}
	&\mathrm{Tr}^kf(t_1,x_1,\dots,t_{m-2k},x_{m-2k})
	\\&\deq \int_{([0,1]\times\R)^{2k}}f(s_1,y_1,\dots,s_{2k},y_{2k},t_1,x_1,\dots, t_{m-2k},x_{m-2k})\notag\\\nonumber
	&\hspace{2em}\times K(s_1-s_2,y_1-y_2)\cdots K(s_{2k-1}-s_{2k},y_{2k-1}-y_{2k})\D {\bbs}\D {\bby}\nonumber.
\end{align}
The following Hu-Meyer's formula (see \cite{hm}, \cite{dk99}, \cite{j06}) connects {\it multiple Stratonovich integrals}  $\mathbb I_m(f)$ with {\it multiple Wiener integrals}  $\I_m(f)$,
\begin{align}\label{e:hm}
	\mathbb I_m(f)=\sum_{k=0}^{[\frac m2]}& \frac{m!}{k!(m-2k)!2^k}\,  \I_{m-2k}( \mathrm {Tr}^k  \hat f \,), ~~ f\in\mathscr H^{ \otimes m},
\end{align}
as long as the right-hand side is well-defined, i.e., $\mathrm{Tr}^k  \hat f\in \mathscr H^{\otimes(m-2k)}$ for $k=0, 1, \dots [m/2]$. In this case, we have
\begin{equation}\label{e:Stra-norm}
	\|f\|^2_{\mathbf S_m} \deq \E[|\mathbb I_m(f)|^2] = \sum_{k=0}^{[\frac m2]} \frac1{(m-2k)!}  \left( \frac{m!}{k!2^k}\right)^2 \|\mathrm {Tr}^k  \hat f \|^2_{\mathscr H^{\otimes(m-2k)}}. 
\end{equation}


To end this section, we introduce the stochastic Fubini theorem which will be used in the proof of Proposition \ref{prop:FK}. Note that stochastic Fubini theorem has been proved in different contexts (see e.g. \cite{DaPrato14,KRT07}), and here we provide a version working for multiple Wiener and Stratonovich integrals.

Let $(\mathbf X, \mathcal X, \mu)$ be a measurable space with $\mu(\mathbf X)<\infty$. For a measurable function $f:\mathbf X\to \R$, we take the notation \[\mu(f) \deq \int_{\mathbf X} f(x) \mu(\D x).\]  The following result holds for general Gaussian spaces associated with a \emph{separable} Hilbert space.
\begin{prop} \label{prop:Fubini}
	For $m\in \NN$,  let  $h: \mathbf X\to \mathscr H^{\otimes m}$ be a measurable mapping such that
	\[\mu(\E[|\I_m(h)|^2])=\mu(\|h\|^2_{\mathscr H^{\otimes m}})<\infty. \]Then the following stochastic Fubini theorem holds,
	\begin{equation}\label{e:Fubini}
		\mu(\mathbf I_m (h))= \mathbf I_m(\mu(h)).
	\end{equation}
	Similarly, if we assume 
	\[\mu\left(\E[|\mathbb I_m(h)|^2]\right)=\mu(\|h\|^2_{\mathbf S_m})=\sum_{k=0}^{[\frac m2]} \frac1{(m-2k)!}  \left( \frac{m!}{k!2^k}\right)^2 \mu\left( \|\mathrm {Tr}^k  \hat h \|^2_{\mathscr H^{\otimes(m-2k)}}\right)<\infty,\]
	we also have 
	\begin{equation}\label{e:Fubini'}
		\mu(\mathbb I_m (h))= \mathbb I_m(\mu(h)).
	\end{equation}
	
\end{prop}
\proof  Let $\{e_i, i\in\mathbb N\}$ be an orthonormal basis of $\mathscr H$. Suppose \[h(x) =\sum_{i_1,\dots, i_m=1}^\infty \alpha_{i_1,\dots, i_m} (x) e_{i_1}\otimes \cdots
\otimes e_{i_m}.\] Then the condition $\mu(\|h\|^2_{\mathscr H^{\otimes m}})<\infty$ implies  $\sum_{i_1, \dots, i_m}\mu( \alpha_{i_1, \dots, i_m}^2)<\infty$. This further implies 
\begin{equation}\label{e:fubini}
	\mu(h)=\sum_{i_1, \dots, i_m} \mu(\alpha_{i_1, \dots, i_m}) e_{i_1}\otimes \cdots \otimes e_{i_m},
\end{equation}
where the series on the right-hand side converges in the Hilbert space $\mathscr H^{\otimes m}$ due to Jensen's inequality and the assumption $\mu(\mathbf X)<\infty$.

Consider $g\in\mathscr H^{\otimes m}$ and suppose 
\[g=\sum_{i_1, \dots, i_m=1}^\infty \beta_{i_1, \dots, i_m} e_{i_1}\otimes \cdots \otimes e_{i_m}\] with $ \|g\|^2_{\mathscr H^{\otimes m}} =\sum_{i_1, \dots, i_m} \beta^2_{i_1, \dots, i_m} <\infty$. We have
\begin{align}\label{e:test-fun}
	\mu(\langle h, g\rangle_{\mathscr H^{\otimes m}})&=\mu\left(\sum_{i_1,\dots,i_m} \alpha_{i_1,\dots,i_m} \beta_{i_1,\dots,i_m}\right)
	\\&	=  \sum_{i_1,\dots,i_m} \mu(\alpha_{i_1,\dots,i_m}) \beta_{i_1,\dots,i_m}=  \langle \mu(h), g\rangle_{\mathscr H^{\otimes m}}\nonumber, 
\end{align}
where the second equality follows from classical Fubini's theorem, noting that \[\mu\left(\sum_{i_1,\dots,i_m} |\alpha_{i_1,\dots,i_m}\beta_{i_1,\dots,i_m}|\right) \le \mu\left(\left(\sum_{i_1,\dots,i_m} \alpha_{i_1,\dots,i_m}^2\right)^{1/2} \right)\left(\sum_{i_1,\dots,i_m}\beta_{i_1,\dots,i_m}^2 \right)^{1/2}<\infty,\]
and the last equality follows from \eqref{e:fubini}. 

By \eqref{e:test-fun} we have $\E[\mu(\I_m(h)) \I_m(g)] = \E[\I_m(\mu(h)) \I_m(g)]$ for any $g\in \mathscr H^{\otimes m}$, and thus the desired equality \eqref{e:Fubini} holds.  Finally, equation \eqref{e:Fubini'} follows from \eqref{e:Fubini} and\eqref{e:hm}. \qed

\subsection {Mild Stratonovich solution} \label{sec:she}

In this subsection, we consider the fractional SHE  \eqref{SHE1} with initial value $u(0,x)=u_0(x), x\in{\mathbb R}$. We call $\{u(t,x)\}_{(t,x)\in \R_+\times \R}$ a {\it mild Stratonovich  solution} of  \eqref{SHE1} if $u(t,x)$ is an adapted process (with respect to the filtration $(\mathcal F_t)_{t\ge0}$ where $\mathcal F_t=\sigma(W(s,x), 0\le s\le t, x\in \R)$) such that $\E[|u(t,x)|^2]<\infty$ for all $(t,x)\in\R_+\times \R$ and satisfies the following integral equation 
\begin{equation}\label{solution-1}
	u(t,x)=\int_{\mathbb R}\gg(t,x-y)u_0(y)\D y+\be\int_0^t\int_{\mathbb R}\gg(t-r,x-y)u(r,y) W(\D r,\D y),
\end{equation}
where $\gg(t,x)$ is the density function of the $\rho$-stable process $\{X_t, t\ge 0\}$ with $X_0=0$ and the integral on the right-hand side  is understood in the sense of Stratonovich which will be specified below.  

Assuming $u_0(x)\equiv1$, we have 
\begin{align}\label{solution-2}
	u(t,x)=1+\be\int_0^t\int_{\mathbb R}\gg({t-r},x-y)u(r,y)W(\D r,\D y).
\end{align}
Iterate the equation and formally we have a series expansion for a mild Stratonovich solution:
\begin{equation}\label{chaos-solution}
	\begin{aligned}
		u(t,x)=&1+\sum_{m=1}^\infty\be^m\int_{[0,t]_<^m}\int_{{\mathbb R}^m}
		\prod_{i=1}^{m}\gg(t_{i+1}-t_{i},x_{i+1}-x_{i})W(\D t_i,\D x_i)\\
		=&1+\sum_{m=1}^\infty\be^m\int_{[0,t]^m}\int_{{\mathbb R}^m} \gg_m(\bbt, \bbx;t,x)
		\prod_{i=1}^{m}W(\D t_i,\D x_i)\\
		=&1+\sum_{m=1}^\infty \beta^m \mathbb I_m(\gg_m(\cdot; t,x)),
	\end{aligned}
\end{equation}
where
\begin{align}\label{e:gm}
	&\gg_m(\bbt, \bbx;t,x)\deq\prod_{i=1}^{m}\gg(t_{i+1}-t_{i},x_{i+1}-x_{i})\mathbf 1_{[0,1]^m_<}(t_1, \dots, t_m)
\end{align} 
with \[[0,t]_<^m\deq\{0<t_1<\cdots<t_m<t\} \text{ and } t_{m+1}\deq t, x_{m+1}\deq x,\] 
and on the right-hand side are \emph{multiple Stratonovich integrals}. The series in \eqref{chaos-solution} converges in $L^1(\Omega)$ under the condition \eqref{e:con} (see Proposition \ref{prop:L2-conv-spde}). We also refer to \cite{HNS2011, Song2017} for  an equivalent description.

Knowing that $u(t,x)$ solves a stochastic heat equation, we now aim to derive a Feynman-Kac type representation for $u(t,x)$.  Given a path of the $\rho$-stable process $X$,  the following Wiener integral
\begin{equation}\label{e:cal-I}
	\mathcal I_{t,x}^{\ep }\deq\ \int_{0}^{t}\int_{\mathbb{R}} p_\ep(X_{t-r}^x-y) W(\D r, \D y),
\end{equation}
is well defined, where $X_{s}^x\deq X_{s}+x$ and $
p_\ep(x) \deq \frac{1}{\sqrt{2\pi \ep}} \e^{-\frac{x^2}{2\ep}}$ is the heat kernel. 

As for the discrete model, we use $\E_X$ and $\E_W$ to denote the expectations in the probability space of  $X$ and $W$, respectively, and we abuse the notation $\E=\E_W\times \E_X$ in this section. 

\begin{prop}
	\label{prop:appr}
	Assume the condition \eqref{e:con}. Then, for each $\ep>0$, $p_\ep(X_{t-\cdot}^x-\cdot)$ belongs to $\mathscr {H}$ a.s.\ and the family of random variables $\mathcal I_{t,x}^{\ep}$ defined by (\ref{e:cal-I}) converges in $L^{2}$  to a
	limit denoted by
	\begin{equation}\label{e:def-Itx}
		\mathcal I_{t,x}\deq \int_{0}^{t}\int_{\mathbb{R}}\delta (X_{t-r}^{x}-y)W(\D r , \D y),
	\end{equation}
	where $\delta(X_{t-\cdot}^{x}-\cdot)$ is an $\mathscr H$-valued random variable given by the $L^2$-limit of $p_\ep(X_{t-\cdot}^x-\cdot)$. Moreover, $\delta(X_{t-\cdot}^{x}-\cdot) \in L^2(\Omega_X, \mathcal G, P; \mathscr H)$,  where $(\Omega_X,\mathcal G, P)$ is the probability space of $X$. 
	
	Conditional on $X$, $\mathcal I_{t,x}$ is a Gaussian random variable with mean $0$
	and variance%
	\begin{equation} \label{e:var-I}
		\mathrm{Var}[\mathcal I_{t,x}|X] =\int_{0}^{t}\int_{0}^{t}|r-s|^{2H-2}\delta(X_r-X_s)\D r \D s\,.
	\end{equation}
\end{prop}

\proof 
For $\varepsilon, \sigma>0$, by \eqref{e:inner} we have for $t\le 1$,
\begin{align}
	&\left\langle p_\ep(X_{t-\cdot}^x-\cdot)\I_{[0,t]} , p_{\si}(X_{t-\cdot}^x-\cdot) \I_{[0,t]} \right\rangle _{\mathscr{H}} \notag\\
	=&\int_0^t \int_0^t \int_{\R} |s-r|^{2H-2} p_{\ep
	}(X_{s}^{x}-y)p_{\sigma}(X_{r}^{x}-y) \D y \D s\D r\label{e4}\\
	=& \int_0^t\int_0^t |s-r|^{2H-2}p_{\ep+\si}(X_s-X_r) \D s \D r\notag. 
\end{align}
By \eqref{e:g-scaling}, we have  \[\E[p_{\ep+\si}(X_s-X_r)] =\int_{\R} p_{\ep+\sigma}(y) \gg(|r-s|,y) \D y \lesssim |r-s|^{-1/\rho}. \]
Thus, the condition \eqref{e:con} yields
\[\E \left\langle p_\ep(X_{t-\cdot}^x-\cdot) , p_{\si}(X_{t-\cdot}^x-\cdot) \right\rangle _{\mathscr{H}}\lesssim  \int_0^t\int_0^t |s-r|^{2H-2-\frac1\rho} \D s \D r<\infty,\]
hence  $p_\ep(X_{t-\cdot}^x-\cdot)$ belongs to $\mathscr{H}$ for all $\ep>0$ almost surely and 
\begin{align}\label{e:2.13}
	\E[\mathcal I_{t,x}^{\ep}\mathcal I_{t,x}^{\sigma}]	&=\E\left\langle p_\ep(X_{t-\cdot}^x-\cdot), p_\si(X_{t-\cdot}^x-\cdot)\right\rangle _{\mathscr{H}}\\&=\E \int_0^t\int_0^t |s-r|^{2H-2}p_{\ep+\si}(X_s-X_r) \D s \D r\nonumber.
\end{align}%
As $(\varepsilon,\sigma)\to 0$,  by the dominated convergence theorem we have\begin{align*}\lim_{(\varepsilon,\sigma)\to 0}  \E[\mathcal I_{t,x}^{\ep}\mathcal I_{t,x}^{\sigma}]&=\E \int_0^t\int_0^t |s-r|^{2H-2}\delta (X_s-X_r) \D s\D r\\&= C \int_0^t\int_0^t |s-r|^{2H-2-\frac1\rho} \D s \D r<\infty
\end{align*}for some proper constant $C$. As a consequence, we have
$\lim_{(\varepsilon,\sigma)\to 0} \E \Big[\Big( \mathcal I_{t,x}^{{\varepsilon}}-\mathcal I_{t,x}^{\sigma}\Big) ^{2}\Big] = 0$.

This implies that,  as $\ep$ converges to zero, $p_{\ep}(X_{t-\cdot}^x-\cdot)$ is a Cauchy sequence in $L^2(\Omega_X,\mathcal G, P; \mathscr H)$, and the limit is denoted by $\delta(X_{t-\cdot}^x-\cdot)$. Therefore,  $\mathcal I_{t, x}^\ep$ converges to $\mathcal I_{t,
	x}$ given in \eqref{e:def-Itx}  in $L^2$. Finally,  (\ref{e:var-I}) can be proven by a similar argument.
\qed

\begin{remark}
	By the computation above, it is clear that we have
	\[\E[ |\mathcal I_{t,x}|^2 ] = C \int_0^t \int_0^t |s-r|^{2H-2-\frac1\rho}\D s\D r,  \]
	which is finite iff \eqref{e:con} holds. It turns out that \eqref{e:con} is also sufficient (and necessary, of course) for $\E[\e^{\mathcal I_{t,x}}]<\infty$ by \eqref{e:var-I} and Proposition \ref{prop:exp-int} below. 
\end{remark}

\begin{prop}\label{prop:exp-int} Under the condition \eqref{e:con}, we have 
	\begin{equation}\label{e:exp-int}
		\E\left[\exp\left(\beta \int_0^1\int_0^1 |r-s|^{2H-2}\delta(X_r-X_s)\D r\D s \right)\right] <\infty \text{ for all } \beta>0.
	\end{equation}
\end{prop}
\proof The proof essentially follows from that of \cite[Theorem 3.3]{Song2017} with some minor modifications. Taylor expansion and Fubini's theorem yield
\begin{equation} \label{e:taylor}
	\begin{aligned}
		&\E\left[\exp\left(\beta \int_0^1\int_0^1 |r-s|^{2H-2}\delta(X_r-X_s)\D r\D s \right)\right]\\
		&=\sum_{m=0}^\infty \frac{1}{m!} \left(\frac\beta{2\pi}\right)^m \int_{[0,1]^{2m}}\int_{\R^m} \prod_{i=1}^m |s_{2i}-s_{2i-1}|^{2H-2} \E\left[\e^{\i\sum_{i=1}^m \xi_i (X_{s_{2i}}-X_{s_{2i-1}})}\right]\D\bbxi \D \bbs\\
		&=\sum_{m=0}^\infty \frac{1}{m!} \left(\frac\beta{2\pi}\right)^m \sum_{\sigma\in \mathcal P_{2m}}\int_{[0,1]_<^{2m}} \int_{\R^m}\prod_{i=1}^m |s_{\sigma(2i)}-s_{\sigma(2i-1)}|^{2H-2}\\&\hspace{15em}\times \E\left[\e^{\i\sum_{i=1}^m \xi_i (X_{s_{\si(2i)}}-X_{s_{\si(2i-1)}})}\right]\D\bbxi \D \bbs,
	\end{aligned}
\end{equation}
where we recall that $[0,1]_<^{2m}=\{0<s_1<s_2<\cdots < s_{2m}<1\}$ and $\mathcal P_{2m}$ is the set of all permutations on $\llb 2m\rrb$.    Here we use the fact  $\delta(x)=\frac1{2\pi}\int_{\R} \e^{\i\xi x} d\xi$.   We remark that  the  integral appearing in \eqref{e:exp-int} is defined as follows $$\int_0^1\int_0^1 |r-s|^{2H-2}\delta(X_r-X_s)\D r\D s:= \lim_{\varepsilon\to 0}\int_0^1\int_0^1 |r-s|^{2H-2}p_\varepsilon(X_r-X_s)\D r\D s,$$
where the limit is taken in $L^2$,  and the computations in \eqref{e:taylor} can be made rigorous by a limiting argument (see, e.g., the proof of \cite[Theorem 4.1]{Song2017}).

Let $\sigma\in \mathcal P_{2m}$ be arbitrarily chosen and fixed. Consider $(s_1, \dots, s_{2m})\in [0,1]^{2m}_<$, and  for each pair $(s_{\sigma(2i-1)}, s_{\sigma(2i)})$ ,  we let $t^*_{2i}\deq s_{\sigma(2i-1)}\vee s_{\sigma(2i)}$ and  $t^*_{2i-1}$ be the unique $s_j$ which is the closest point to $t^*_{2i}$ from the left. Then clearly, we have \[\prod_{i=1}^m |s_{\sigma(2i)}-s_{\sigma(2i-1)}|^{2H-2}\le  \prod_{i=1}^m |t^*_{2i}-t^*_{2i-1}|^{2H-2}.\]

Noting that \[\sum_{i=1}^m \xi_i (X_{s_{\si(2i)}}-X_{s_{\si(2i-1)}})=\sum_{j=2}^{2m} \eta_j (X_{s_j}-X_{s_{j-1}})\] where each $\eta_j$ is a linear combination of $\xi_i$'s for  $i=1,\dots, m.$  Then we have 
\begin{align*}
	\E\left[\e^{\i\sum_{i=1}^m \xi_i (X_{s_{\si(2i)}}-X_{s_{\si(2i-1)}})}\right]= \prod_{j=2}^{2m} \e^{-c_\rho |\eta_j|^{\rho} (s_j-s_{j-1})}\le \prod_{i=1}^m \e^{-c_\rho |\tilde \eta_i|^{\rho} (t_{2i}^*-t_{2i-1}^*)},
\end{align*}
where the inequality holds because we only keep the factors resulting from the characteristic function of $X_{t_{2i}^*}-X_{t_{2i-1}^*}$ for $i\in\llb m\rrb$ and drop all the others. Here, for each $i\in\llb m\rrb$, we define $\tilde \eta_i=\eta_j$ where $j$ is the unique index such that $s_j=t^*_{2i}=s_{\sigma(2i)}\vee s_{\sigma(2i-1)}$. 

Thus, we have 
\begin{equation}\label{e:estimation-moment}
	\begin{aligned}
		&\int_{[0,1]_<^{2m}} \int_{\R^m}\prod_{i=1}^m |s_{\sigma(2i)}-s_{\sigma(2i-1)}|^{2H-2} \E\left[\e^{\i\sum_{i=1}^m \xi_i (X_{s_{\si(2i)}}-X_{s_{\si(2i-1)}})}\right]\D\bbxi \D \bbs\\
		&\le \int_{[0,1]_<^{2m}} \int_{\R^m} \prod_{i=1}^m |t^*_{2i}-t^*_{2i-1}|^{2H-2}  \prod_{i=1}^{m} \e^{-c_\rho |\tilde \eta_i|^{\rho} (t^*_{2i}-t^*_{2i-1})} \D\bbxi \D \bbs\\
		&= C^m \int_{[0,1]_<^{2m}} \prod_{i=1}^m |t^*_{2i}-t^*_{2i-1}|^{2H-2-\frac1\rho} \D \bbs\\
		&\le   \frac{C^m }{\Gamma\left(m(2H-\frac1\rho)+1\right)}.
	\end{aligned}
\end{equation}
where the equality follows from  the fact that the Jacobian determinant $| [\partial \xi_i/\partial{\tilde \eta_{j}}]_{m\times m}|$ is one and that $\int_{\R} \e^{-a |\xi|^\rho t} \D \xi=C a^{-1/\rho}$ for $a>0$,  and the last step follows from Lemma \ref{lem:ineq-gamma}.

Therefore, by \eqref{e:taylor} and \eqref{e:estimation-moment}, and noting that $|\mathcal P_{2m}|=(2m)!$,  we have 
\begin{align*}&\E\left[\exp\left(\beta \int_0^1\int_0^1 |r-s|^{2H-2}\delta(X_r-X_s)\D r\D s \right)\right]\\&\le \sum_{m=0}^\infty C^m\beta^m \frac{(2m)!}{m! \Gamma\left(m(2H-\frac1\rho)+1\right)} <\infty,
\end{align*}
where the last step follows from Stirling's formula and $2H-\frac1\rho>1$. 
\qed

\begin{remark}
	In the proof of \cite[Theorem 3.3]{Song2017}, the Markov property instead of the independent increment property of $X$ was invoked, and as a consequence, \cite[Lemma 2.2]{Song2017} was needed therein. As shown in the proof of Proposition \ref{prop:exp-int}, if we  utilise the independent increment property of $X$,  \cite[Lemma 2.2]{Song2017} is no longer needed, and this is important for the proof of uniform $L^1$-bound for rescaled partition function $\ZZ_\om^{(N)}(\hat \be_N, k)$ in Section \ref{sec:bd-polymer}.    
\end{remark}

\begin{remark}
	Note that when $\rho=2$, the condition $\theta>\frac12$ in \eqref{e:con}  becomes $H>\frac34$. Thus, Proposition \ref{prop:exp-int} is consistent with \cite[Theorem 6.2]{HNS2011}. 
\end{remark}


\begin{prop}\label{prop:FK} Under the condition \eqref{e:con}, the mild Stratonovich solution  given in \eqref{chaos-solution} has the following Feynman-Kac representation:
	\begin{equation}\label{e:F-K}
		u(t, x) = \E_X\left[ \exp\left(\beta \int_0^t\int_{{\mathbb R}} \delta(X_{t-r}^{x}-y)
		W(\D r,\D y)\right)\right].  
	\end{equation}
\end{prop} 
\proof  First, we prove the integrability of the right-hand side of \eqref{e:F-K}. For all $\beta\in \R$,
\begin{equation}\label{e:exp-int'}
	\begin{aligned}
		& \E_W
		\left[\E_X\left[ \exp\left(\beta \int_0^t\int_{{\mathbb R}} \delta(X_{t-r}^{x}- y)
		W(\D r,\D  y)\right)\right]  \right]\\
		&=\E_X\left[\E_W\left[ \exp\left(\beta \int_0^t\int_{{\mathbb R}} \delta(X_{t-r}^{x}- y)
		W(\D r,\D  y)\right)\right]  \right]\\
		&=
		\E_X\left[ \exp\left(\frac 12 \beta^2 t^{2H-\frac1\rho} \int_0^1 \int_0^1 |r-s|^{2H-2}\delta(X_r-X_s)\D r\D s \right)\right] \\
		&<\infty,
	\end{aligned}
\end{equation}
where the second equality follows from \eqref{e:var-I} and the self-similarity of $X$, and the last step follows from Proposition \ref{prop:exp-int}.  Thus $u(t, x)$ given by \eqref{e:F-K} is well-defined and is $L^p$-integrable for all $p>0$. 

The Feynman-Kac formula \eqref{e:F-K}  now follows from \eqref{chaos-solution} and the following equation:
\begin{equation}\label{e:F-K-chaos}
	\begin{aligned}
		&\E_X\left[ \left(\int_0^t\int_{{\mathbb R}} \delta(X_{t-r}^{x}- y)
		W(\D r,\D  y)\right)^m\right]\\
		&= \int_{[0,t]^m}\int_{\R^{m}} \E_X \left[ \prod_{i=1}^m \delta(X_{t-r_i}^
		x- y_i) \right] W(\D r_1, \D  y_1)\dots W(\D r_m, \D  y_m)  \\
		&= m! \int_{[0,t]^m}\int_{{\mathbb R}^{m}} \gg_m( \bbr, \bby;t,x)
		\prod_{i=1}^{m}W(\D r_i,\D y_i),
	\end{aligned}
\end{equation}
where the first equality follows from \eqref{e:II-m} and the stochastic Fubini theorem (Proposition \ref{prop:Fubini}). \qed
\begin{remark} \label{rem:g-m-tr} A direct corollary  of  Propositions  \ref{prop:exp-int} and \ref{prop:FK} is 
	\begin{equation}
		\|\gg_m\|_{\mathbf S_m}<\infty,
	\end{equation}
	where $\gg_m$ is given in \eqref{e:gm} and $\|\cdot\|_{\mathbf S_m}$   in \eqref{e:Stra-norm}. Indeed, 
	\begin{align*}
		\|\gg_m\|^2_{\mathbf S_m}=\E[|\mathbb I_m(\gg_m)|^2]=\left(\frac1{m!}\right)^2 \E_W\left[ \Big(\E_X[\mathcal I_{t,x}^m]\Big)^2\right]\le \left(\frac1{m!}\right)^2  \E[\mathcal I_{t,x}^{2m}]<\infty,
	\end{align*}
	where $\mathcal I_{t,x}$ is given in \eqref{e:def-Itx}, the second equality follows from \eqref{e:F-K-chaos}, and the last inequality holds since $\E[\e^{\beta\mathcal I_{t,x}}]<\infty$ for all $\beta\in\R$ by \eqref{e:exp-int'}.  Thus, $\mathrm{Tr}^k \hat\gg_m\in \mathscr H^{\otimes(m-2k)}$ by \eqref{e:Stra-norm}. Moreover, by \eqref{e:F-K-chaos} we have
	\[\sum_{m=1}^\infty \E[|\mathbb I_m(\gg_m)|]\le \sum_{m=0}^\infty \frac{1}{m!} \E[\left|\mathcal I_{t,x}\right|^{m}] =\E[{\e^{\left|\mathcal I_{t,x}\right|}}]\le 2\E[\e^{\mathcal I_{t,x}}]<\infty.\]
	This implies that $\sum_{m=1}^\infty \mathbb I_m(\gg_m)$ converges in $L^1$, which is stated in the following proposition.
\end{remark}

\begin{prop}\label{prop:L2-conv-spde}   Under the condition \eqref{e:con}, the mild Stratonovich  solution of \eqref{solution-2} given in \eqref{chaos-solution}  converges in $L^1$. 
\end{prop}

\subsection{Mild It\^o-Skorohod solution}\label{sec:she1}

Mild It\^o-Skorohod solution is defined similarly as for mild Stratonovich solution in Section \ref{sec:she}, and the only difference is that the stochastic integral in \eqref{solution-2} is now understood in the Skorohod sense. More precisely, still assuming $\tilde u_0(x)\equiv1$, the It\^o-Skorohod solution $\tilde u(t,x)$ of \eqref{SHE1} satisfies
\begin{align}\label{solution-3}
	\tilde  u(t,x)=1+\be\int_0^t\int_{\mathbb R}\gg({t-r},x-y)\tilde  u(r,y) \diamond W(\D r,\D y),
\end{align}
and its chaos expansion is given by
\begin{equation}\label{chaos-solution1}
	\begin{aligned}
		\tilde  u(t,x)
		=&1+\sum_{m=1}^\infty\be^m\int_{[0,t]^m}\int_{{\mathbb R}^m} \gg_m(\bbt, \bbx;t,x)
		W(\D t_1,\D x_1)\diamond \cdots \diamond W(\D t_m,\D x_m)\\
		=&1+ \sum_{m=1}^\infty \beta^m \I_m(\gg_m(\cdot; t,x)),
	\end{aligned}
\end{equation}
where $\gg_m$ is given in \eqref{e:gm} and the stochastic integrals above are in the Skorohod sense. A sufficient and necessary condition for the existence and uniqueness of the It\^o-Skorohod solution is (see, e.g., \cite{hn09,Song2017}), 
\begin{equation}\label{e:chaos-sum-norm}
	\sum_{m=1}^\infty \beta^m m! \left\|\hat \gg_m(\cdot;t,x)\right\|^2_{\mathscr H^{\otimes m}}<\infty,
\end{equation}
where $\hat \gg_m$ is the symmetrization of $\gg_m$:
\begin{equation}\label{e:g-hat}
	\hat \gg_m(\bbt, \bbx;t,x)=\frac1{m!} \sum_{\sigma\in \mathcal P_m} \gg_m (\bbt_\sigma, \bbx_\sigma; t,x),
\end{equation}
with $\bbt_\sigma\deq  (t_{\sigma(1)}, \dots, t_{\sigma(m)})$ and $\bbx_\sigma\deq (x_{\sigma(1)}, \dots, x_{\sigma(m)})$. 

By \cite[Theorem 5.3]{Song2017}, there exists a unique mild It\^o-Skorohod solution for \eqref{SHE1} on $\R^d$ if 
\begin{equation}\label{e:con-sko}
	\int_{\R^d} \frac{1}{1+|\xi|^\rho}\D \xi <\infty \Longleftrightarrow \rho>d.
\end{equation}

\begin{remark}\label{rem:sko-sol}
	Note that the requirement $\rho \le 2$ together with \eqref{e:con-sko} forces us to consider the case $d=1$ only. By \eqref{e:con-sko}, we also need to assume $\rho>1$ in Theorem \ref{thm:main1}. 
\end{remark}

We shall provide an alternative proof for \eqref{e:chaos-sum-norm} under the condition $H\in(1/2, 1], \rho \in (1,2]$, which will be easier to be adapted to estimate the moments of the terms arising from the  partition function $\tilde Z_\om^{(N)}$ (see Remark \ref{rem:bd-discrete} below).

\begin{prop}\label{prop:g-m}
	Assume $H\in(1/2, 1]$ and $\theta=H-\frac1{2\rho}>0$. Then there exists a constant $C$ such that for all $(t,x)\in [0,1]\times \R$, 
	\begin{equation}\label{e:m-chaos}
		\|\hat \gg_m(\cdot; t,x)\|_{\mathscr H^{\otimes m}} \le  C^m (m!)^{H-1}\left( \frac{\Gamma(\theta/H)^{m}}{ \Gamma(m \theta/H+1 )} \right )^{H},
	\end{equation}
	where $\hat \gg_m$ is given in \eqref{e:g-hat}. Moreover, \eqref{e:chaos-sum-norm} holds under the condition \eqref{e:con'}. 
\end{prop}

\proof 
Recalling \eqref{e:inner}, we have 
\begin{align*}
	\|\hat \gg_m(\cdot;t,x)\|^2_{\mathscr H^{\otimes m}} &=\int_{[0,1]^{2m}}\int_{\R^m} \prod_{i=1}^m |r_i-s_i|^{2H-2}\hat\gg_m(\bbr, \bbx;t,x)\hat \gg_m(\bbs, \bbx;t,x) \D \bbx \D \bbr\D\bbs\\
	&\le C_H^m \int_{\R^m} \left(\int_{[0,1]^m} |\hat\gg_m(\bbs, \bbx;t,x)|^{\frac1H}\D \bbs\right)^{2H}\D\bbx\\
	&\le C_H^m \left(\int_{[0,1]^m} \left(\int_{\R^m} |\hat\gg_m(\bbs, \bbx;t,x)|^{2} \D\bbx\right)^{\frac1{2H}}\D\bbs\right)^{2H}
\end{align*}
where the first inequality follows from Lemma \ref{lem:HLS} and the second one follows from the Minkowski inequality.  Recalling the definition \eqref{e:g-hat} of $\hat \gg_m$ and the definition \eqref{e:gm} of $\gg_m$, we have
\begin{align*}
	|\hat\gg_m(\bbs, \bbx;t,x)|^{2}&=\frac1{(m!)^2}\left|\sum_{\sigma\in\mathcal P_m}\gg_m(\bbs_\sigma, \bbx_\sigma;t,x)\right|^2 \\
	&= \frac1{(m!)^2}\left|\sum_{\sigma\in\mathcal P_m}\prod_{i=1}^{m}\gg(s_{\sigma(i+1)}-s_{\sigma(i)},x_{\sigma(i+1)}-x_{\sigma(i)})\mathbf 1_{[0,1]^m_<}(s_{\sigma(1)}, \dots, s_{\sigma(m)})\right|^2\\
	&= \frac1{(m!)^2}\sum_{\sigma\in\mathcal P_m}\prod_{i=1}^{m} \left|\gg(s_{\sigma(i+1)}-s_{\sigma(i)},x_{\sigma(i+1)}-x_{\sigma(i)})\right|^2 \mathbf 1_{[0,1]^m_<}(\bbs_\sigma),
\end{align*}
where we use the convention $s_{\sigma(m+1)}=s, x_{\sigma(m+1)}=x$, and hence
\begin{align*}
	&\|\hat \gg_m(\cdot;t,x)\|^2_{\mathscr H^{\otimes m}} \\&\le C_H^m \left(\int_{[0,1]^m} (m!)^{-\frac{1}{H}}\left(\int_{\R^m} \sum_{\sigma\in\mathcal P_m}\prod_{i=1}^{m} \left|\gg(s_{\sigma(i+1)}-s_{\sigma(i)},x_{\sigma(i+1)}-x_{\sigma(i)})\right|^2 \mathbf 1_{[0,1]^m_<}(\bbs_\sigma)\D\bbx\right)^{\frac1{2H}}\D\bbs\right)^{2H} \\
	&=C_H^m \left(\int_{[0,1]_<^m} (m!)^{1-\frac{1}{H}}\left(\int_{\R^m} \prod_{i=1}^{m}|\gg(s_{i+1}-s_{i},x_{i+1}-x_{i})|^2\D\bbx\right)^{\frac1{2H}}\D\bbs\right)^{2H},
\end{align*}
Noting that by \eqref{e:g-bd}, we have that for all $t>0$, 
\begin{equation*}\label{e:g2-int}
	\int_{\R} |\gg(t,x)|^2 \D x \le C t^{-\frac1\rho},
\end{equation*}
where $C$ is a positive constant depending only on $\rho$. Combining the last two estimates together with  Lemma \ref{lem:ineq-gamma}, we get 
the desired inequality \eqref{e:m-chaos}. 

Finally, by \eqref{e:m-chaos} and Stirling's formula we get that \eqref{e:chaos-sum-norm} holds if $\rho>1$. 
\qed 

\begin{remark}\label{rem:-rho} The scaling limit of long-range random walk in i.i.d. random environment was considered in \cite{CAR-SUN-ZYG}. In this situation,  the condition $\rho>1$ is also a necessary condition for \eqref{e:chaos-sum-norm}. Indeed, now we have
	\begin{align*}
		\|\hat \gg_m(\cdot;t,x)\|^2_{\mathscr H^{\otimes m}} &=\int_{[0,1]^{m}}\int_{\R^m} |\hat \gg_m(\bbs, \bbx;t,x)|^2 \D \bbx \D\bbs\\
		&=\frac{1}{m!} \int_{[0,1]^m_<} \int_{\R^m}  \prod_{i=1}^{m}|\gg(s_{i+1}-s_{i},x_{i+1}-x_{i})|^2\D\bbx \D\bbs. 
	\end{align*}
	By  the scaling property \eqref{e:g-scaling}, we have 
	\[\int_{\R} |\gg(s,x)|^2 \D x  = A s^{-\frac1\rho},\]
	where $A=\int_\R |\gg(1, x)|^2 \D x$ is finite by  \eqref{e:g-bd}.  Thus, we get
	\begin{align*}
		\|\hat \gg_m(\cdot;t,x)\|^2_{\mathscr H^{\otimes m}} = \frac{A^m }{m!}\int_{[0,1]_<^m} \prod_{i=1}^m (s_{i+1}-s_i)^{-\frac1\rho}\D \bbs,
	\end{align*}
	which is finite only if $\rho>1$. One can also check that \eqref{e:chaos-sum-norm} holds by Stirling's formula if $\rho>1$.
\end{remark}

\begin{remark}\label{rem:bd-discrete}
	From the proof of Proposition \ref{prop:g-m}, we can see that the estimation \eqref{e:m-chaos} still holds if we replace $\gg(t,x)$ by its upper bound  $C(t\|x\|^{-\rho-1}) \wedge t^{-1/\rho}$ (see \eqref{e:g-bd}). This fact together with the  upper bound for $P_n(k)$ given in \eqref{e:Pn-bd''}   will be used to estimate the moments of ${\mathbf S}_m^{(N)}$ for the Skorohod case in Section \ref{sec:bd-polymer}.
\end{remark}

\begin{remark}\label{rem:con-two-solutions} Assume $H\in(\frac12,1]$ and $\rho\in(0,2]$. For \eqref{SHE1}, by the analysis in this section (or by Theorem 4.6 and Theorem 5.3 in \cite{Song2017}), we know that the condition for the existence of a Stratonovich solution is $\theta =H-\frac1{2\rho}>\frac12$ and the condition for the It\^o-Skorohod case is $\rho>1.$ Note that for the Stratonovich case,  $\theta>\frac12$ yields $\rho>\frac1{2H-1}\ge1$ and moreover  $\theta>\frac12,\rho\le 2$ implies $H\ge\frac34$. Thus, it requires more restrictive condition for the existence of a Stratonovich solution. This is because the $L^2$-norm of $\mathbb I_m(\gg_m)$  is strictly bigger than $\mathbf I_m(\gg_m)$  due to the extra trace terms  appearing in $\|\gg_m\|^2_{\mathbf S_m}$ (see \eqref{e:Stra-norm}). 
	
\end{remark}

\section{On weak convergences}\label{sec:proof}
In this section,  we aim to prove the main results Theorems \ref{thm:main} and \ref{thm:main1}.  Recall that $\theta=H-\frac1{2\rho}$. For $N\in\NN$, under the scaling $\beta\rightarrow\hat\beta_N=\beta N^{-\theta}$,  the rescaled  partition function $Z_\om^{(N)}$  given by \eqref{partition2} is 
\begin{align}\label{sptp}
	Z_\om^{(N)}(\hat\beta_N,N^{1/\rho}x_0)=\E_S\left[\exp\left(\hat\beta_N \sum_{n=1}^N \om\left(n, S^{(N+1,k)}_n\right)\right)\right],
\end{align}
where  $x_0\in\RR$ is fixed such that  $S^{(N+1,k)}_{N+1}=N^{1/\rho}x_0=k\in \ZZ$.  To simply the notation, we use $S$ to denote the backward random walk $S^{(N+1,k)}$. Taylor expansion yields 
\begin{equation}\label{esptp}
	\begin{aligned}
		Z_\om^{(N)}(\hat\beta_N,N^{1/\rho}x_0)&=\sum_{m=0}^\infty\frac{1}{m!}{\mathbb S}_m^{(N)},
	\end{aligned}
\end{equation}
where  ${\mathbb S}_m^{(N)}$ is the $m$th moment of the Hamiltonian (with a factor $\hat \beta_N$) given by 
\begin{equation}\label{e:S-m}
	\begin{aligned}
		{\mathbb S}_m^{(N)}=&\E_S\left[\left(\hat\beta_N \sum_{n=1}^N \om(n,S_n)\right)^m \right]=\hat\beta_N^m\E_S\left[\sum_{\bbn\in \llb N \rrb^m}\om(n_1,S_{n_1})\cdots\om(n_m,S_{n_m})\right].
	\end{aligned}
\end{equation}

For every $\bbn=(n_1, \dots, n_m)\in \llb N \rrb^m$, the components can be arranged in an increasing order and the resulted  sequence is denoted by  $\bbn^*=(n_1^*, \cdots, n_m^*)$ with $n^*_1\le n^*_2\le\cdots\le n^*_m$.  For each $\bbn\in \llb N\rrb^m$, there is a permutation $\sigma$ of $\llb N \rrb$ such that $n_i^*=n_{\sigma_i}$ for $i=1, \dots, m$. Denote 
\begin{equation}\label{e:P*}
	\begin{aligned} 
		&P_{\bbn^*}
		=P_{\bbn^*}(N^{1/\rho}x_0; k_1, \dots,  k_m)\\
		&=P_{n^*_2-n^*_1}( k_{\sigma_2}- k_{\sigma_1}) \cdots P_{n_m^*-n_{m-1}^*}( k_{\sigma_m}- k_{\sigma_{m-1}}) P_{(N+1)-n_m^*}(N^{1/\rho}x_0-k_{\sigma_m}) ,
	\end{aligned}
\end{equation}
where we use the convention $P_0(0)=1$ and $P_0( k)=0$ for $ k\neq 0$. We remark that $P_{\bbn^*}$ is \emph{symmetric} in its $m$ arguments.  The $m$th moment in \eqref{e:S-m} can be written as
\begin{align}\label{4-1}
	{{\mathbb S}}_m^{(N)}=\hat\beta^m_N\sum_{n_1, \dots, n_m\in \llb N \rrb}\sum_{ k_1,\dots, k_m\in\ZZ}\om(n_1, k_1)\cdots\om(n_m, k_m)P_{\bbn^*}(N^{1/\rho} x_0; k_1,\dots, k_m),
\end{align}
recalling that $S$ is  a backward random walk with $S_{N+1}= N^{1/\rho}x_0$.

A similar calculation can be done for the  partition function $\tilde \ZZ_\om^{(N)}$ given in \eqref{partition3}. Note that the exponential is actually a Wick exponential of $\beta\sum_{i=1}^N \om(i, S_i)$ conditional on the random walk $S$. Hence, applying \eqref{e:wick-exp}, we get the following resemblance of \eqref{esptp} and \eqref{e:S-m} respectively,
\begin{equation}\label{esptp'}
	\begin{aligned}
		\tilde Z_\om^{(N)}(\hat\beta_N,N^{1/\rho}x_0)&=\sum_{m=0}^\infty\frac{1}{m!}{\mathbf S}_m^{(N)},
	\end{aligned}
\end{equation}
where  
\begin{equation}\label{e:S-m'}
	\begin{aligned}
		{\mathbf S}_m^{(N)}=&\E_S\left[\left(\hat\beta_N \sum_{n=1}^N \om(n,S_n)\right)^{\diamond m} \right]\\
		=&\hat\beta_N^m\E_S\left[\sum_{\bbn\in \llb N \rrb^m}:\!\om(n_1,S_{n_1})\cdots\om(n_m,S_{n_m})\!:\right]\\
		=&  \beta^m N^{-m(\theta+\frac1\rho)}\sum_{n_1, \dots, n_m\in \llb N \rrb}\sum_{ k_1,\dots, k_m\in\ZZ}:\!\om(n_1, k_1)\cdots\om(n_m, k_m)\!:\\
		&\qquad \qquad \qquad \qquad\qquad \qquad  \times \left(N^{m/\rho}P_{\bbn^*}(N^{1/\rho} x_0; k_1,\dots, k_m)\right).
	\end{aligned}
\end{equation}

Now, we introduce the following general U-statistics  $\mathbf{I}_m^{(N)}(f)$  which will be proven to converge weakly to a multiple Wiener integral.  For  $f\in \mathcal B^{\otimes m}$ where we recall that $\mathcal B$ is given in \eqref{eq-B},  we denote
\begin{align}\label{U-stat-I}  
	\mathbf{I}^{(N)}_m(f)&\deq N^{-m(\theta+\frac1\rho)} \sum_{n_1, \dots, n_m\in \llb N \rrb}\sum_{ k_1,\dots, k_m\in\ZZ}\bigg[:\!\om(n_1, k_1)\cdots\om(n_m, k_m)\!:\notag\\
	&\qquad \qquad \qquad \qquad\qquad \qquad \qquad \qquad  \times {\mathcal A_N(f)}\left(\frac{n_1}{N}, \frac{ k_1}{N^{1/\rho}},\cdots, \frac{n_m}{N},\frac{ k_m}{N^{1/\rho}}\right)\bigg]\\
	&\deq N^{-m(\theta+\frac1\rho)} \sum_{\bbn\in\llb N\rrb^m}\sum_{ \bbk\in \ZZ^m}:\!\om_{\llb m\rrb}\!: \mathcal A_N(f)(\bbt,\bbx),\notag
\end{align}
where $:\!\om_{\llb m\rrb}\!:\deq :\!\om_1\cdots \om_m\!:$ is the physical Wick product of $\om_i\deq\om(n_i,  k_i)$, ${\mathcal A_N(f)}$ is given in \eqref{con-exp}, and 
\begin{equation}\label{e:nk-tx}
	\begin{aligned}
		&\bbt=(t_1, \dots, t_m)\deq\bbn/N=(n_1/N,\cdots, n_m/N), \\
		&\bbx=(x_1, \dots, x_m)\deq\bbk/N^{\frac1\rho}=( k_1/N^{\frac1\rho}, \cdots,  k_m/N^{\frac1\rho}). 
	\end{aligned}
\end{equation}

In order to write $\mathbf S_m^{(N)}$ defined in \eqref{e:S-m'} in the form of $\mathbf I_m^{(N)}(f)$ given in \eqref{U-stat-I}, it suffices to extend the domain of the following function of $(\bbt,\bbx)$
\[N^{m/\rho}P_{N\bbt }(N^{1/\rho} x_0; N^{1/\rho}x_1,\dots, N^{1/\rho} x_m):=N^{m/\rho}P_{\bbn^*}(N^{1/\rho} x_0; k_1,\dots, k_m)\] to the whole $[0,1]_<^m\times\R^m$ in a natural way, i.e., we define
 \begin{equation}\label{e:tildeP}
		\tilde P_m( \bbt, \bbx)\deq N^{m/\rho}P_{\bbn^*}(N^{1/\rho} x_0; k_1,\dots, k_m), \text{ if } (\bbt, \bbx)\in \Big[\frac{{\bbn^*}}{N},\frac{{\bbn^*}+\overrightarrow{1}_m}{N}\Big)\times\Big[\frac{{ \bbk}}{N^{1/\rho}},\frac{{\bbk}+\overrightarrow{1}_m}{N^{1/\rho}}\Big),
\end{equation} 
In this way, we can write $\mathbf S_m^{(N)} =  \beta^m \mathbf I_m^{(N)} (\tilde P_m)$.

 The rest of this Section is organized as follows. In Section \ref{sec:U-stat}, we prove the joint weak convergence for $U$-statistics $\I_m^{(N)} (f)$  defined by \eqref{U-stat-I}. The joint weak convergence of $m$th moment ${\mathbb S}_m^{(N)}$ (resp. $\mathbf S_m^{(N)}$) appearing in the rescaled partition function $Z_\om^{(N)}(\hat \beta_N, N^{1/\rho}x_0)$  (resp.  $(\tilde Z_\om^{(N)}(\hat \beta_N, N^{1/\rho}x_0)$) is established in Section \ref{sec:wc-partition}. The $L^p$-bounds of $Z_\om^{(N)}(\hat \be_N, k)$ and $\tilde Z_\om^{(N)}(\hat \be_N, k)$ are obtained in Section \ref{sec:bd-polymer}.
\subsection{Weak convergence of U-statistics}\label{sec:U-stat}

We first introduce the following lemma which provides a uniform (in $N$) bound for the second moment of  $\I_m^{(N)} (f)$  with $f\in \mathcal B^{\otimes m}$, which will be used to prove  Propositions \ref{thm-lim-1} and \ref{prop:m-integral-I}.

\begin{lemma}\label{lem:m-mom-bound}  
	Assume the disorder $\om$ is given as in Section \ref{sec:model}. Then, there exists a constant $C$ such that for all  $f\in \mathcal B^{\otimes m}$ and $m, N\in\mathbb N$,
	\begin{equation}\label{e:I-1-bound}
		\E\left[\left(\I_m^{(N)}(f)\right)^2\right]\le C^m m! \|\hat f\|^2_{\mathcal B^{\otimes m}},
	\end{equation}
	where  $\hat f$ is the symmetrization of $f$.
\end{lemma}

\proof  For $m=1$, we have   
\begin{align}\label{I-1}
	\begin{split}
		&\E\left[\left(\I_1^{(N)}(f)\right)^2\right]\\
		=&N^{-2(\theta+1/\rho)}\sum_{n, n'\in \llb N \rrb} \sum_{ k, k'\in \ZZ} \E[\om(n, k)\om(n', k')]{\mathcal A_N(f)}\left(\frac{n}{N},\frac{ k}{N^{1/\rho}}\right){\mathcal A_N(f)}\left(\frac{n'}{N},\frac{ k'}{N^{1/\rho}}\right)\\
		=& N^{-2(\theta+1/\rho)}\sum_{n,n'\in \llb N \rrb} \sum_{ k\in \ZZ} \gamma(n-n'){\mathcal A_N(f)}\left(\frac{n}{N},\frac{ k}{N^{1/\rho}}\right){\mathcal A_N(f)}\left(\frac{n'}{N},\frac{k}{N^{1/\rho}}\right).
	\end{split}
\end{align}

For $t\in\R$, let 
\begin{equation}\label{e:gamma-N}
	\gamma_N(t)\deq N^{2-2H}\gamma\left(\left[ |t|N\right]\right).
\end{equation} Then,  by \eqref{e:gamma-bd} and \eqref{e:gamma-lim}, we have for all $t\in\R$
\begin{equation}\label{e:gamma-N-bd}
	\gamma_N(t)\lesssim |t|^{2H-2} \text{ and } \lim\limits_{N\to\infty}\gamma_N(t)=|t|^{2H-2}.
\end{equation}
Then, using the notations $t=n/N, t'=n'/N, x= k/N^{1/\rho}$, we have
\begin{equation}\label{e:I-N-f}
	\begin{aligned}
		\E\left[\left(\I_1^{(N)}(f)\right)^2\right]&= N^{-2-1/\rho} \sum_{n,n'\in\llb N \rrb}\sum_{k\in\mathbb Z} N^{2-2H} \gamma\left(n-n'\right) {\mathcal A_N(f)}\left(t,  x\right) {\mathcal A_N(f)}\left(t',  x\right)  \\
		&\lesssim \int_0^1\int_0^1\int_{\RR} |t-t'|^{2H-2} \left|{\mathcal A_N(f)}(t, x){\mathcal A_N(f)}(t',{x})\right| \D x \D t\D t'\\
		&= \|{\mathcal A_N(f)}\|^2_{\mathcal B}\lesssim \|{f}\|^2_{\mathcal B},
	\end{aligned}
\end{equation}
where the last inequality follows from Lemma \ref{lem:Jensen}.

Recall that $\theta=H-\frac1{2\rho}$ given in \eqref{e:con}.  For general $m$, 
\begin{align}\label{e:I-m-2-mom}
	\E\left[\left(\I_m^{(N)}(f) \right)^2 \right] =N^{-2m(\theta+1/\rho)} \sum_{\bbn, \bbn'}\sum_{ \bbk,  \bbk'} \E\left[:\!\om_{\llb m\rrb}\!::\!\om'_{\llb m\rrb}\!:\right]\mathcal A_N(f)(\bbt, \bbx)\mathcal A_N(f)(\bbt', \bbx'),
\end{align}
where $:\!\om_{\llb m\rrb}\!:=:\!\prod_{i=1}^m \om(n_i,  k_i)\!:$ and similarly for $:\!\om'_{\llb m\rrb}\!:$. 

Recalling that we have assumed that $\om$ is Gaussian, then by  \eqref{e:E-wick-prod'} we have 
\[\E\left[:\!\om_{\llb m\rrb}\!::\!\om'_{\llb m\rrb}\!:\right] = \E\left[:\!\prod_{i=1}^m \om(n_i,  k_i)\!::\!\prod_{i=1}^m \om(n_i',  k_i')\!:\right] = \sum_{\sigma\in \mathcal P_m} \prod_{i=1}^m\gamma(n_i-n_{\sigma(i)}')  \delta_{k_i k'_{\sigma(i)}}, \]
where the summation $\sum_{\sigma\in \mathcal P_m}$ is taken over the set $\mathcal P_m$ of all the permutations of $\llb m\rrb$. Note that by the symmetry of the summation $\sum_{\bbn, \bbn'}\sum_{ \bbk,  \bbk'}$, we have $\I_m^{(N)}(f)=\I_m^{(N)}(\hat f)$.  Thus, by the symmetry of $\mathcal A_N(\hat f)$, we have
\[\E\left[\left(\I_m^{(N)}(\hat f) \right)^2 \right] =m!\,N^{-2m(\theta+1/\rho)}  \sum_{\bbn, \bbn'}\sum_{ \bbk} \left(\prod_{i=1}^m\gamma(n_i-n_i') \right)\mathcal A_N(\hat f)(\bbt, \bbx)\mathcal A_N(\hat f)(\bbt', \bbx).\]
Then, similarly to the case $m=1$, by Lemma \ref{lem:Jensen} we can get  $\E\left[\left(\I_m^{(N)}(f)\right)^2\right]\le   C^m \|\hat f\|^2_{\mathcal B^{\otimes m}}.$ \qed

Now, we are ready to prove the weak convergence for $\I_m^{(N)}(f)$ defined in \eqref{U-stat-I}.  We first prove the weak convergence for $m=1$. 

\begin{prop}\label{thm-lim-1}
	Let $f: [0,1]\times\RR \to \RR$ be a function in $\mathcal B$. Then, $\I_1^{(N)}(f)$ converges weakly to the Wiener integral $\I_1(f)$ as $N\to\infty$.  Moreover, for any $k\in\NN$ and $f_1, \dots, f_k\in \mathcal B$, we have the joint convergence in distribution:
	\begin{equation}\label{e:joint-con-1}
		\left(\I_1^{(N)}(f_1), \dots, \I_1^{(N)}(f_k)\right)  \wc \left(\I_1(f_1), \dots, \I_1(f_k)\right), \text{ as } N\to\infty. 
	\end{equation}
\end{prop}

\proof First we consider an indicator function $f={\mathbf 1}_{[s,t]\times[ y, z]}$ with $0\le s\leq t \le 1$ and $ y\leq  z$. In this case,  
\begin{align*}
	\I_1^{(N)}(f)=N^{-\theta-1/\rho}\sum_{n\in[Ns,Nt]\cap \NN}~\sum_{ k\in[N^{1/\rho} y,N^{1/\rho} z]\cap \ZZ}\om(n, k).
\end{align*} 
Clearly it has mean zero and  the covariance is
\begin{align*}
	\E\left[\left(\I_1^{(N)}(f)\right)^2\right]
	&=  N^{-2-1/\rho} \sum_{n,n'\in[Ns, Nt]\cap\NN}\sum_{ k\in[N^{1/\rho} y,N^{1/\rho} z]\cap \ZZ} \gamma_N(r-r'),
\end{align*}
where  $\gamma_N$ is given by \eqref{e:gamma-N} and $r=n/N,r'=n'/N$.  By \eqref{e:gamma-N-bd} and the dominated convergence theorem, we get 
\begin{equation}\label{e:I-1f}
	\begin{aligned}
		\lim_{N\to\infty} \E\left[\left(\I_1^{(N)}(f)\right)^2\right]&=  \int_s^t\int_s^t \int_y^z  |r-r'|^{2H-2} \D  x\D r\D r'  =\|f\|^2_{\mathscr H}
	\end{aligned}
\end{equation}
where we recall that the inner product  $\la \cdot, \cdot \ra_{\mathscr H}$ is given   in \eqref{e:inner}.  Thus, we have $\I_1^{(N)}(f)\wc \I_1(f)$ as $N\to\infty$, noting that $\I_1^{(N)}(f), N\in \NN$ and  $\I_1(f)$ are Gaussian random variables and that $\I_1(f)$ has zero mean and variance $\|f\|_{\mathscr H}^2$.

Similarly, one can show that \eqref{e:I-1f} holds for each simple function $f$, i.e., $f$ is a linear combination of indicator functions.  For a  general function $f\in\mathcal B$, there exists an approximation sequence of simple functions $\{f^{(n)}\}_{n\in \NN}$ in $\mathcal B$. On the one hand, $\I_1^{(N)}(f^{(n)})$ converges weakly to  $\I_1(f^{(n)})$ as $N\to \infty$   for each fixed $n$.  On the other hand, by Lemma \ref{lem:m-mom-bound}, we have $\I_1^{(N)}(f^{(n)})\to \I_1^{(N)}(f)$
in $L^2$ uniformly in $N$ as $n\to\infty$. Combining the above convergences together with the obvious  $L^2$-convergence of  $\I_1(f^{(n)})\to \I_1(f)$ as $n\to \infty$, by Lemma~\ref{lem:weak-con} we obtain the weak convergence $\I_1^{(N)}(f)\wc \I_1(f)$ as $N\to \infty$ for $f\in\mathcal B$.     The argument can be presented in a diagram:
\[\xymatrixcolsep{10pc} \xymatrix{ \I_1^{(N)}(f^{(n)}) \ar[d]_{d}^{N \to\infty} \ar[r]^{\mathrm{in\ } L^2,\mathrm{ \ uniformly\ in\ } N }_{n \to\infty}& \I_1^{(N)}(f) \ar[d]_{d}^{N\to\infty}\\ \I_1(f^{(n)}) \ar[r]_{n\to\infty}^{\mathrm{in \ } L^2} & \I_1(f))}
\]

Finally, the joint weak convergence \eqref{e:joint-con-1} follows from the linearity of $\I_1^{(N)}$ and the Cram\'er-Wold theorem (see Theorem \ref{thm:cw}). 
\qed

The main result in this section is presented below. 

\begin{prop}\label{prop:m-integral-I}
	For each $m\in \NN$ and $f\in\mathcal B^{\otimes m}$, we have
	\begin{align}\label{e:con-m}
		{\mathbf I}_m^{(N)}(f)\wc  \I_m(f), \text{ as } N\to \infty, 
	\end{align}
	where ${\mathbf I}_m^{(N)}(f)$ is given in \eqref{U-stat-I} and $\I_m(f)$  is an $m$th multiple Wiener integral.  Moreover, for any $k\in\NN$ and $f_1, \dots, f_k$ with $f_i\in \mathcal B^{\otimes l_i}$ with $l_i\in\NN$, we have the joint convergence in distribution:
	\begin{equation}\label{e:joint-con-m}
		\left(\I_{l_1}^{(N)}(f_1), \dots, \I_{l_k}^{(N)}(f_k)\right) \wc \Big(\I_{l_1}(f_1), \dots, \I_{l_k}(f_k)\Big), \text{ as } N\to \infty. 
	\end{equation}
\end{prop}

\proof    We first consider $\I_m^{(N)}(f)$ for $f=h^{\otimes m}$ with $h\in\mathcal B$.  For such $f$, we have
\begin{align*}
	{\mathbf I}_m^{(N)}(f)=&N^{-m(\theta+1/\rho)} \sum_{n_1, \dots, n_m\in\llb N\rrb}\sum_{ k_1, \dots,  k_m \in\ZZ} :\!\om_{\llb m\rrb}\!: \mathcal A_N(h)^{\otimes m}	= :\!\Big({\mathbf I}_1^{(N)}(h)\Big)^{m}\!:\,,
\end{align*}
where the term on the right-hand side is a physical Wick product (see Section \ref{sec:wick}). Note that $\I_1^{(N)}(h)$ is a centered Gaussian random variable with variance bounded by $\|h\|^2_{\mathcal B}$ up to a multiplicative constant by  Lemma \ref{lem:m-mom-bound}. Therefore, for any $q\in\NN$, by the Gaussianity of $\om$ and \eqref{e:I-N-f}, we have 
$$\E\left[\left(\I_1^{(N)}(h)\right)^{2q}\right]\lesssim (2q-1)!!\|h\|^{2q}_{\mathcal B}<\infty.$$
Then we can apply Lemma \ref{lem:weak-con'} together with Proposition \ref{thm-lim-1} and get that,  as $N\to\infty$,
\begin{align*}
	{\mathbf I}_m^{(N)}(f)= :\!\Big({\mathbf I}_1^{(N)}(h)\Big)^m\!:\wc &:\!\Big({\mathbf I}_1(h)\Big)^{m}\!:=\I_m(h^{\otimes m})=\I_m(f).
\end{align*}

To prove the result, it suffices to consider symmetric functions in $\mathcal B^{\otimes m}$. For any symmetric function $f\in \mathcal B^{\otimes m}$,  one can find a sequence $\{f^{(n)}\}_{n\in\NN}$  which are the linear combinations of the functions of the form  $h^{\otimes m}$ such that $f^{(n)}\to f$ in $\mathcal B^{\otimes m}$. By the preceding argument, for each $f^{(n)}$, we have $\I_m^{(N)}(f^{(n)}) \wc \I_m(f^{(n)})$ as $N\to\infty$. Then similarly as in the proof of Proposition \ref{thm-lim-1},  we have, recalling Lemma \ref{lem:m-mom-bound},
\[\xymatrixcolsep{10pc} \xymatrix{ \I_m^{(N)}(f^{(n)}) \ar[d]_{d}^{N \to\infty} \ar[r]^{ \mathrm{in\ } L^2,\mathrm{ \ uniformly\ in\ } N}_{n \to\infty}& \I_m^{(N)}(f) \\ \I_m(f^{(n)}) \ar[r]_{n\to\infty}^{\mathrm{ in \ } L^2} & \I_m(f))}
\]
This together with Lemma~\ref{lem:weak-con} yields $\I_m^{(N)}(f)\wc \I_m(f)$ as $N\to \infty$ for symmetric $f\in\mathcal B^{\otimes m}$.

Now we prove \eqref{e:joint-con-m}. Similar to the preceding step, we first consider the case $f_i= h_i^{\otimes l_i}$ with $h_i\in\mathcal B, i=1, \dots, k$, then Proposition \ref{thm-lim-1} and Lemma \ref{lem:weak-con''}  yield that for all $(a_1,\dots, a_k)\in\R^k$, 	$ \sum_{i=1}^k a_i \I_{l_i}^{(N)}(f_i) \wc \sum_{i=1}^k a_i\I_{l_i}(f_i)$ as $N\to \infty$, which implies \eqref{e:joint-con-m}. The general case follows from a limiting argument together with Lemma~\ref{lem:weak-con}. \qed

\subsection{Weak convergence of rescaled partition functions} \label{sec:wc-partition}

 In this section, we prove  Theorem \ref{thm:main} for the Stratonovich case  and Theorem \ref{thm:main1} for the It\^o-Skorohod case. As explained in Section \ref{sec:main-result}, we shall focus on the proof for the Stratonovich case which also  encompasses the essentials for the It\^o-Skorohod case.

Recall that the so-called $m$th moment ${\mathbb S}_m^{(N)}$ appearing in the rescaled partition function $Z_\om^{(N)}(\hat \beta_N, N^{1/\rho}x_0)$ has the expression~\eqref{4-1}. We first prove the joint weak convergence of moments. 

\begin{prop}\label{m-thm}
	Assume condition \eqref{e:con}. For each $m\in \NN$, we have
	\begin{align}\label{e:con-mom}
		\frac1{m!}{\mathbb S}_m^{(N)}\wc   \beta^m \mathbb I_m(\gg_m), \text{ as } N\to \infty, 
	\end{align}
	where $\gg_m=\gg_m(\bbt, \bbx;1,x_0)$ is given in \eqref{e:gm} and $\mathbb I_m(\gg_m)$ is an $m$th multiple Stratonovich integral.  Moreover, for any $k\in\NN$ and $l_1, \dots, l_k\in\NN$, we have the joint convergence in distribution:
	\begin{equation}\label{e:joint-con-mom}
		\left(\frac{1}{l_1!}\mathbb S_{l_1}^{(N)}, \dots, \frac{1}{l_k!}\mathbb S_{l_k}^{(N)}\right)  \wc \Big(\beta^{l_1}\mathbb I_{l_1}(\gg_{l_1}), \dots, \beta^{l_k}\mathbb I_{l_k}(\gg_{l_k})\Big), \text{ as } N\to \infty. 
	\end{equation}
\end{prop}
\proof We first prove \eqref{e:con-mom}. Note that $\mathrm{Tr}^k \hat\gg_m\in \mathcal H^{\otimes(m-2k)}$ (see \eqref{e:trace-k} for the definition of $\mathrm{Tr}^k \hat\gg_m$) for $k=0,1,\dots, [\frac m2]$ (see Remark \ref{rem:g-m-tr}). Thus, by Proposition \ref{prop:m-integral-I}, we have (we assume $m$ is an odd integer and the analysis for the case of an even number  $m$  is the same),
\begin{equation}\label{e:con-multi}
	\begin{aligned}
		&\bigg(\I_{1}^{(N)}\left(\mathrm{Tr}^{[\frac m2]}\hat\gg_{m}\right),\I_{3}^{(N)}\left(\mathrm{Tr}^{[\frac m2]-1}\hat\gg_{m}\right), \dots, \I_{m}^{(N)}\left(\hat\gg_m\right)\bigg)\\
		&\wc \bigg(\I_{1}\left(\mathrm{Tr}^{[\frac m2]}\hat\gg_{m}\right),\I_{3}\left(\mathrm{Tr}^{[\frac m2]-1}\hat\gg_{m}\right), \dots, \I_{m}\left(\hat\gg_m\right)\bigg), \text{ as } N\to\infty.
	\end{aligned}
\end{equation}

Recalling  \eqref{4-1}, by the symmetry of the summation, we can write ${\mathbb S}_m^{(N)}$ as 
\[{\mathbb S}_m^{(N)} =\hat\beta^m_N\sum_{n_1, \dots, n_m\in \llb N \rrb}\sum_{ k_1,\dots, k_m\in \ZZ}\om_{\llb m \rrb} P_{\bbn^*}, 
\]
where $\hat\beta_N=\beta N^{-\theta},  P_{\bbn^*}$ is given in \eqref{e:P*} and $\om_{\llb m\rrb}=\prod_{i=1}^m \om_i$ with $\om_i=\om(n_i,  k_i)$. It follows from  \eqref{wick} and \eqref{e:expectation'} that
$$ \om_{\llb m \rrb}=\sum_{j=0}^{[\frac m2]}\sum_{\substack{B\subset {\llb m \rrb}\\ |B|=m-2j}} :\!\om_
B\!: \E\left[\om_{{\llb m \rrb}\backslash B}\right].$$
Therefore, 
\begin{equation}\label{eq-m-chaos}
	\begin{aligned}
		{\mathbb S}_m^{(N)}=&\hat\beta_N^m\sum_{n_1, \dots, n_m} \sum_{ k_1,\dots,  k_m} \Bigg(\sum_{j=0}^{[\frac m2]}\sum_{\substack{B\subset {\llb m \rrb}\\ |B|=m-2j}} :\!\om_B\!: \E\left[\om_{{\llb m \rrb}\backslash B}\right] \Bigg)  P_{\bbn^*}\\
		=& \sum_{j=0}^{[\frac m2]}\sum_{\substack{B\subset {\llb m \rrb}\\ |B|=m-2j}}\Bigg(\hat\beta_N^m\sum_{n_1, \dots, n_m} \sum_{ k_1,\dots,  k_m}  :\!\om_B\!: \E\left[\om_{{\llb m \rrb}\backslash B}\right]   P_{\bbn^*}\Bigg).  
	\end{aligned}
\end{equation}

For each $j=0,1, 
\dots, [\frac m2]$, all the terms in the summation $\sum_{\substack{B\subset {\llb m \rrb}\\ |B|=m-2j}}$ are equal  by the symmetry.  We fix $j\in\{0, 1, \dots, [\frac m2]\}$. Note that there are in total ${m \choose m-2j}$ subsets in $\llb m \rrb$ with cardinality  $m-2j$. Without loss of generality, we assume $B=\{2j+1,\dots,m\}$. Then by \eqref{e:expectation'}, we have
\begin{align}\label{part-ex}
	\E[\om_{\llb m \rrb\backslash B}]=\E[\om_{\llb 2j\rrb}] = \sum_{V} \prod_{\{\ell_1,\ell_2\}\in V} \E[\om_{\ell_1}\om_{\ell_2}],
\end{align}
where the sum $\sum_V$ is taken over all pair partitions of $\llb2j\rrb$. By the symmetry again,  the summations \[\sum_{n_1, \dots, n_m} \sum_{k_1,\dots,k_m}  :\!\om_B\!:  \prod_{\{\ell_1,\ell_2\}\in V} \E[\om_{\ell_1}\om_{\ell_2}]    P_{\bbn^*}\]
coincide with each other for all the partitions $V$, and note that 
there are in total $(2j-1)!!$ pair partitions of $\llb2j\rrb$.  Thus, we have the following equation:
\begin{align}\label{e:e:hm'}
	\mathbb S_m^{(N)}& =\sum_{j=0}^{[\frac m2]} {m\choose{m-2j}} (2j-1)!! ~\left[\hat\beta_N^m\sum_{n_1, \dots, n_m} \sum_{ k_1,\dots,  k_m}  :\!\om_B\!:  \prod_{\ell=1}^{j} \E[\om_{2\ell-1}\om_{2\ell}]    P_{\bbn^*}\right]\notag\\
	& =\sum_{j=0}^{[\frac m2]} \frac{m!}{j!(m-2j)! 2^j} ~\left[\hat\beta_N^m \sum_{n_1, \dots, n_m} \sum_{ k_1,\dots,  k_m}  :\!\om_B\!:  \prod_{\ell=1}^{j} \E[\om_{2\ell-1}\om_{2\ell}]    P_{\bbn^*}\right],
\end{align}
where we recall $B=\{2j+1, \dots, m\}$.

Compare \eqref{e:e:hm'} with Hu-Meyer's formula \eqref{e:hm} and also note \eqref{e:con-multi}. It then follows easily that  in order to prove \eqref{e:con-mom},  it suffices to show, setting  $\beta=1$  without loss of generality and hence $\hat \beta_N=N^{-\theta}$, 
\begin{equation}\label{e:AN}
	\begin{aligned}
		Y_N&\deq \frac{1}{m!} N^{-m\theta} \sum_{n_1, \dots, n_m} \sum_{ k_1,\dots,  k_m}  :\!\om_B\!:  \prod_{\ell=1}^{j} \E[\om_{2\ell-1}\om_{2\ell}]    P_{\bbn^*} -  \I_{m-2j}^{(N)}(\mathrm{Tr}^{j} \hat\gg_m)\\
		&= \frac{1}{m!} N^{-(m-2j) (\theta+\frac1\rho)} \sum_{n_{2j+1}, \dots, n_m} \sum_{ k_{2j+1},\dots,  k_m}  :\!\om_B\!: \\
		&\qquad\times\bigg( N^{\frac m\rho}N^{-2j(\theta+\frac1\rho)}\sum_{n_i: i\in \llb 2j \rrb} \sum_{ k_i: i\in \llb 2j\rrb} \prod_{\ell=1}^{j} \E[\om_{2\ell-1}\om_{2\ell}]    P_{\bbn^*}\bigg) -  \I_{m-2j}^{(N)}(\mathrm{Tr}^{j} \hat\gg_m)\\
		&=  \I_{m-2j}^{(N)}\bigg(N^{-2j(\theta+\frac1\rho)}\sum_{n_i: i\in \llb 2j \rrb} \sum_{ k_i: i\in \llb 2j\rrb}\prod_{\ell=1}^{j} \E[\om_{2\ell-1}\om_{2\ell}]   \left( \tilde P_m/m!\right)\bigg)   - \I_{m-2j}^{(N)}(\mathrm{Tr}^{j} \hat\gg_m)\\
		&\to  0 \text{ in $L^2$ as } N \to \infty,
	\end{aligned}
\end{equation}
for all $j=0, 1, \dots, [\frac m2]$, where we recall that  $\I_{m-2j}^{(N)}(f)$ is defined in \eqref{U-stat-I} and $\tilde P_m$ is given in \eqref{e:tildeP}.  In the rest of the proof, we may abuse the notation $N^{m/\rho}  P_{\bbn^*}$ for  $\tilde P_m(\bbt,\bbx)$ where it is appropriate.

We prove \eqref{e:AN} for $j=0$, and the general case can be proved in a similar spirit.  When $j=0$, we have
 \[Y_N= \I_m^{(N)}\left( \tilde P_{m}/m! - \hat \gg_m \right). \]
Denoting $\mathfrak D_N=\mathfrak D_N(\bbn, \bbk) = \tilde P_{m}/m! - \mathcal A_N(\hat \gg_m)(\bbn/N, \bbk/N^{1/\rho}), \mathfrak D'_N= \mathfrak D_N(\bbn', \bbk')$ and $\tilde{\mathfrak D}'_N= \mathfrak D_N(\bbn', \bbk)$, we have
\begin{equation}\label{e:AN'}
	\begin{aligned}
		\E[Y_N^2]=&N^{-2m(\theta+\frac1\rho)} \sum_{\bbn, \bbn' \in \llb N\rrb^m}\sum_{\bbk, \bbk'\in \ZZ^m}\E\left[:\!\om_{\llb m \rrb}\!::\!\om'_{\llb m \rrb}\!:\right] \mathfrak D_N\mathfrak D'_N\\
		=& m! N^{-2m(\theta+\frac1\rho)} \sum_{\bbn, \bbn' \in \llb N\rrb^m}\sum_{\bbk\in \ZZ^m}\prod_{i=1}^m \gamma(n_i-n_i') \mathfrak D_N\tilde{\mathfrak D}'_N\\
		\lesssim & m! \int_{[0,1]^{2m}\times \R^m} \prod_{i=1}^m |t_i-t_i'|^{2H-2} \mathfrak D_N\tilde{\mathfrak D}'_N \D \bbt \D \bbt' \D \bbx
	\end{aligned}
\end{equation}
where the second equality follows from \eqref{e:E-wick-prod'} and we use the notations that link the discrete summations and continuum integrals:  $t_i=n_i/N, t'_i=n'_i/N$ and   $x_i= k_i/N^{1/\rho}$. In light of Lemma~\ref{lem:Jensen}, we have
\begin{equation}\label{e:est-YN}
	\E[Y_N^2] \lesssim m! \int_{[0,1]^{2m}\times \R^m} \prod_{i=1}^m |t_i-t_i'|^{2H-2} 
	\left|\mathbf D_N\tilde{\mathbf D}'_N \right|\D \bbt \D \bbt' \D \bbx,  
\end{equation}
where $\mathbf D_N=\mathbf D_N(\bbn, \bbk) = \tilde P_{m}/m! - \hat \gg_m(\bbt, \bbx)$ and $\tilde{\mathbf D}_N'=\mathbf D_N(\bbn', \bbk).$

Inspired by \cite{CAR-SUN-ZYG}, we decompose $I\deq [0,1]^{2m}\times \R^m$ as  $I\cap(I_1\cup I_2 \cup I_3)$, where
\[I_1\deq \bigcap_{i=1}^m \left(\Big\{|t_i-t'_{i}|>\ep \Big\}\cap \{|x_i|<M\}\right),\]
\[I_2\deq \bigcup_{i=1}^m \{|t_i-t'_{i}|\le \ep\}, \text{ and } I_3\deq \bigcup_{i=1}^m\{|x_i|\ge M\},\]
for some fixed $\ep, M>0$. On $I\cap I_1$, $\prod_{i=1}^m|t_i-t_i'|^{2H-2} |\mathbf D_N\tilde{\mathbf D}'_N|$ is uniformly bounded and  converges to 0 by the local limit theorem, and hence the integral converges to 0 by the dominated convergence theorem. Then, we can get $\lim_{N\to\infty}\E[Y_N^2]=0$, if we can show the integral on $I\cap I_2$ (resp. $I\cap I_3$) can be arbitrarily small if we choose $\ep$ sufficiently small (resp. $M$ sufficiently large). 

By \eqref{e:Pn-bd} and \eqref{e:g-scaling}, we have 
\begin{equation}\label{e:est-gg}
	\begin{aligned}
		& \left(N^{m/\rho} P_{\bbn^*}\right)\left(N^{m/\rho} P_{(\bbn')^*}\right) \lesssim \Big((t^*_2-t^*_1)\cdots(t_{m}^*-t_{m-1}^*)( 1-t_m^*) \Big)^{-1/\rho} N^{m/\rho} P_{(\bbn')^*},\\
		&\hat\gg_m(\bbt, \bbx;1,x_0)\hat\gg_m(\bbt', \bbx;1,x_0) \lesssim \Big((t^*_2-t^*_1)\cdots(t_{m}^*-t_{m-1}^*)(1-t_m^*)\Big)^{-1/\rho} \hat\gg_m(\bbt',\bbx;1,x_0),
	\end{aligned}
\end{equation}
and similarly for the cross terms $N^{m/\rho} P_{\bbn^*} \hat\gg_m(\bbt',\bbx)$ and $N^{m/\rho} P_{(\bbn')^*} \hat\gg_m(\bbt,\bbx)$, where $t_i^*=t_{\sigma(i)}$ for some permutation $\sigma$ on $\llb m\rrb$ such that $t_1^*\le t_2^*\le \cdots \le t_m^*$.  These inequalities and the fact for all $\bbt\in[0,1]^m$, \[\int_{\R^m} N^{m/\rho} P_{\bbn^*}\D \bbx =  \int_{\R^m} (m!) \hat \gg_m\D \bbx  =1,\] 
yield that $\int_I \prod_{i=1}^m|t_i-t_i'|^{2H-2} |\mathbf D_N\tilde{\mathbf D}'_N| \D \bbt\D \bbt' \D \bbx$ is bounded uniformly in $N$ by (up to a multiplicative constant)
\begin{equation}\label{e:Cm}
	C_m\deq \frac1{m!}  \int_{[0,1]^{2m}} \prod_{i=1}^m |t_i-t_i'|^{2H-2} \Big((t^*_2-t^*_1)\cdots(t_{m}^*-t_{m-1}^*)(1-t_m^*)\Big)^{-1/\rho} \D \bbt \D \bbt', 
\end{equation}
which is finite due to the condition  $2H-2>-1$ and $-1/\rho>-1$.  This implies that the integral  on $I\cap I_2$ can be arbitrarily small by choosing $\ep$ sufficiently small. Similarly,  for the integral on $I\cap I_3$, 
\begin{align*}
	&	\int_{I\cap I_3} \prod_{i=1}^m |t_i-t_i'|^{2H-2}  |\mathbf D_N\tilde{\mathbf D}'_N |\D \bbt\D \bbt' \D \bbx\\&\lesssim  C_m\Bigg(P\left\{\sup_{0\le t\le 1} |X_t|\ge M\Big| X_1= x_0 \right\}+ P\left\{\max_{0\le m\le N} |S_m|\ge N^{1/\rho} M\Big| S_N= N^{1/\rho} x_0 \right\}\Bigg). 
\end{align*}
Noting that \[P\left\{\max_{0\le m\le N} |S_m|\ge N^{1/\rho} M\Big| S_N= N^{1/\rho} x_0 \right\}\to P\left\{\sup_{0\le t\le 1} |X_t|\ge M\Big| X_1= x_0 \right\} \text{ as } N\to\infty,\] we can make the integral on $I\cap I_3$ as small uniformly in $N$ as we want  by choosing $M$ sufficiently large. This proves \eqref{e:AN} for $j=0$. 

For general $j=0, 1,\dots, [\frac m2]$, using the same argument leading to \eqref{e:est-YN}, we have
\begin{equation}\label{e:est-YN'}
	\E[Y_N^2]\lesssim (m-2j)! \int_{[0,1]^{2(m-2j)}\times \R^{m-2j}} \prod_{i=1}^{m-2j} |t_i-t_i'|^{2H-2}\left|\mathbf D_N\tilde{\mathbf D}'_N \right| \D \bbt \D \bbt' \D \bbx,
\end{equation}
where now 
\begin{equation}\label{e:DN}
	\begin{aligned}
		&\mathbf D_N=\mathbf D_N(\bbn, \bbk)\\
		&=N^{-2j\theta}N^{(m-2j)\frac 1\rho}\sum_{n_i: i\in \llb 2j \rrb} \sum_{ k_i: i\in \llb 2j\rrb}\prod_{\ell=1}^{j} \E[\om_{2\ell-1}\om_{2\ell}]   ( P_{\bbn^*}/m!)-\mathrm{Tr}^j\hat \gg_{m}\\
		&\deq \mathcal Q_N- \mathrm{Tr}^j\hat \gg_{m}
	\end{aligned}
\end{equation}
and $\tilde{\mathbf D}_N'=\mathbf D_N(\bbn', \bbk)$.  Then, we estimate the right-hand side of \eqref{e:est-YN'} in a similar way as for the case $j=0$, i.e., we split $[0,1]^{2(m-2j)}\times \R^{m-2j}$ as the union of its restrictions on $I_1, I_2$ and $I_3$ and then analyse the integral restricted on $I_i$'s separately for $i=1, 2, 3$. The analysis for the integral on $I_1$ is the same as for $j=0$.  To argue that the integral on $I_2$ and $I_3$ can be arbitrarily small if we choose $\ep$ sufficiently small and $M$ sufficiently large, it suffices to find a uniform (in N) upper bound for 
\begin{equation}\label{e:term1}
	A\deq \int_{[0,1]^{2(m-2j)}\times \R^{m-2j}} \prod_{i=1}^{m-2j} |t_i-t_i'|^{2H-2} \mathrm{Tr}^j\hat\gg_m(\bbt, \bbx;1,x_0)\mathrm{Tr}^j\hat\gg_m(\bbt', \bbx;1,x_0) \D \bbt \D \bbt' \D \bbx
\end{equation} 
and
\begin{equation}\label{e:term1}
	B\deq \int_{[0,1]^{2(m-2j)}\times \R^{m-2j}} \prod_{i=1}^{m-2j} |t_i-t_i'|^{2H-2} \mathcal Q_N\tilde{\mathcal{Q}}'_N \D \bbt \D \bbt' \D \bbx
\end{equation} 
where $\mathcal Q_N$ is given in \eqref{e:DN} and $ \tilde{\mathcal{Q}}'_N =\mathcal Q_N(\bbn', \bbk).$ This is true,  noting that Remark \ref{rem:g-m-tr} yields
\[A= \frac{1}{(m-2j)!} \E\left[\left|\I_{m-2j}(\mathrm{Tr}^j\hat \gg_m)\right|^2\right] \le \E\left[\left|\mathbb I_{m}( \gg_m)\right|^2\right]<\infty,\]
and similarly, recalling that $\mathbb S_m^{(N)}$ is given in \eqref{4-1}, we have $B\lesssim \E[|\mathbb S_m^{(N)}|^2]$ which has a uniform upper bound by \eqref{e:bd-partition}. 

In the above, we have assumed  $m$ is an odd integer. If $m$ is even, the analysis is the same except for the case $j=\frac m2$. 
In this case, we have
\begin{align*}
	Y_N=N^{-m\theta}\sum_{n_i: i\in \llb m \rrb} \sum_{ k_i: i\in \llb m \rrb}\prod_{\ell=1}^{m/2} \E[\om_{2\ell-1}\om_{2\ell}]   (P_{\bbn^*}/m!)   - \mathrm{Tr}^{m/2} \hat\gg_m,
\end{align*}
which is deterministic and  converges to 0 by the local limit theorem and a similar argument proving $\lim_{N\to\infty}\E[Y_N^2]=0$ for $j=0$. 

Finally, one can prove the weak convergence for linear combinations of $\left(\mathbb S_{l_1}^{(N)}, \dots, \mathbb S_{l_k}^{(N)}\right)$ in a similar way, and hence \eqref{e:joint-con-mom} holds due to Theorem \ref{thm:cw}.  \qed

By the proof  of Proposition \ref{m-thm}, in particular the part proving  \eqref{e:AN} for $j=0$ under the condition $H\in(1/2, 1]$ and $\rho \in (1,2]$, we can get a parallel result  for the It\^o-Skorohod case which is stated below. Recall that $\mathbf S_m^{(N)}$ is given in \eqref{e:S-m'}. 
\begin{prop}\label{m-thm'}
	Assume $H\in(1/2, 1], \rho \in (1,2]$. For each $m\in \NN$, we have
	\begin{align*}\label{e:con-mom'}
		\frac1{m!}{\mathbf S}_m^{(N)}\wc   \beta^m \mathbf I_m(\gg_m), \text{ as } N\to \infty, 
	\end{align*}
	where $\gg_m=\gg_m(\bbt, \bbx;1,x_0)$ is given in \eqref{e:gm} and $\mathbf I_m(\gg_m)$ is an $m$th multiple Wiener integral.  Moreover, for any $k\in\NN$ and $l_1, \dots, l_k\in\NN$, we have the joint convergence in distribution:
	\begin{equation*}\label{e:joint-con-mom'}
		\left(\frac{1}{l_1!}\mathbf S_{l_1}^{(N)}, \dots, \frac{1}{l_k!}\mathbf S_{l_k}^{(N)}\right)  \wc \Big(\beta^{l_1}\mathbf I_{l_1}(\gg_{l_1}), \dots, \beta^{l_k}\mathbf I_{l_k}(\gg_{l_k})\Big), \text{ as } N\to \infty. 
	\end{equation*}
\end{prop}

Now we are ready to prove our main results.

{\it Proof of Theorem \ref{thm:main}.}  By Proposition \ref{m-thm}, we have  for all $M\in \NN$,
\begin{align*}
	Z_\om^{(N,M)}\deq \sum_{m=0}^M\frac{1}{m!}{\mathbb S}_m^{(N)} \wc 
	\mathcal Z^{(M)} \deq \sum_{m=0}^M \beta^m \mathbb I_m ( \gg_m(\cdot;1,x_0)), \text{ as } N\to \infty.
\end{align*}
Recalling Proposition \ref{prop:L2-conv-spde} which yields the $L^1$-convergence of $\mathcal Z^{(M)}$ to $\mathcal Z$ as $M\to \infty$,  we only need to show that $Z_\om^{(N,M)}$ converges to $Z_\om^{(N)}=Z_\om^{(N)}(\hat\beta_N,N^{1/\rho}x_0)$ given in \eqref{esptp} in probability uniformly in $N$ as $M\to\infty$ by Lemma \ref{lem:weak-con}. This follows from \eqref{e:uniform-convergence} obtained in Section \ref{sec:bd-polymer}. 

{\it Proof of Theorem \ref{thm:main1}.}  The proof is the same as that of Theorem \ref{thm:main} except that Proposition~\ref{m-thm} and \eqref{e:uniform-convergence} are replaced by Proposition \ref{m-thm'} and \eqref{e:uniform-convergence'}, respectively.

\subsection{$L^p$-bounds of  rescaled partition functions}\label{sec:bd-polymer}
 In this section, we study the $L^p$-bounds of $Z_\om^{(N)}(\hat \be_N, k)$ and $\tilde Z_\om^{(N)}(\hat \be_N, k)$.
We first deal with the Stratonovich case to find $L^1$-bound of $Z_\om^{(N)}(\hat \be_N, k)$, and   then as a consequence we obtain the $L^2$-bound for the It\^o-Skorohod case.

 Recalling that $\om$ is Gaussian, we have
\begin{align*}
	&\E\left[\exp\left(N^{-\theta} \sum_{n=1}^N \om(n, S_n)\right)\right]\\&=\E\left[ \exp\left(N^{-\theta} \sum_{n=1}^N \sum_{ k\in\ZZ} \om(n,  k) \mathbf 1_{[S_n= k]}\right)\right]\\
	&=\E\left[ \exp\left(\frac12 N^{-2\theta} \sum_{n,n'=1}^N \sum_{ k, k'\in\ZZ} \gamma(n- n')\delta_{kk'} \mathbf 1_{[S_{n}= k]}1_{[S_{n'}= k']}  \right)\right]\\
	&=\E\left[ \exp\left(\frac12 N^{-2\theta} \sum_{n,n'=1}^N \gamma(n-n')\mathbf 1_{[S_{n}=S_{n'}]}  \right)\right].
\end{align*}

In this section, we shall prove a discretised version of Proposition \ref{prop:exp-int}, i.e., to show that the above exponential moment is uniformly bounded.  As  a consequence, we get a discretised version of Proposition \ref{prop:L2-conv-spde} (see eq. \eqref{e:uniform-convergence}).

Noting that $S_n\in \ZZ$ and using the identity 
\[ \mathbf 1_{[S_{n}=S_{n'}]} =\frac{1}{2\pi}\int_{-\pi}^\pi \e^{\i(S_{n}-S_{n'})\lambda} \D \lambda, \]
we have 
\begin{align*}
	&\E\left[ \exp\left(\frac12 N^{-2\theta} \sum_{n,n'=1}^N \gamma(n-n')\mathbf 1_{[S_{n}=S_{n'}]}  \right)\right]\\
	&=\sum_{m=0}^\infty \frac{1}{m!} (4\pi)^{-m}N^{-2m\theta} \sum_{n_1,\dots, n_m=1}^N\sum_{n_1',\dots, n'_m=1}^N \int_{[-\pi, \pi]^m} \prod_{i=1}^m \gamma(n_i-n_i') \E\left[ \e^{\i \sum_{i=1}^m (S_{n_i}-S_{n_i'})\lambda_i }\right]\D \boldsymbol{\lambda}.
\end{align*}

Using changes of variables $n_i=Nt_i, n_i'=Nt_i'$ and $\lambda_i=u_i/N^{1/\rho}$ for $i=1, \dots, m$, we have for each term on the right-hand side of the above equation,
\begin{equation}\label{e:BNm}
	\begin{aligned}
		A_m\deq \frac{1}{m!} (4\pi)^{-m}  N^{-2m} \sum_{\bbn,\bbn'}& \prod_{i=1}^m \left[N^{2-2H}\gamma(N(t_i-t_i'))\right]\\
		&\times  \int_{[- \pi N^{1/\rho} ,  \pi N^{1/\rho} ]^m} \E\left[\e^{ \i \sum_{i=1}^m (S_{n_i} -S_{n_i'})u_i/N^{1/\rho}}\right] \D \bbu .
	\end{aligned}
\end{equation}

Recall that $S_n=S_0+\sum_{j=1}^nY_j$ and that $\psi(u) = \E[\e^{\i u Y_j}]$ is the characteristic function of $Y_j$. 
The 1-lattice distribution of $Y_j$ implies that,  for any $\ep>0$, one can find a positive constant $c$ such that (see \cite[eq. (5.14)]{rosen1990})
\begin{equation}\label{e:bd-psi}
	|\psi(u/N^{1/\rho}) |\le \e^{-c|u|^{\rho-\ep}/N}, \text{ for } 1\le |u| \le \pi N^{1/\rho}.
\end{equation}
Thus,  for $n\in\llb N\rrb$, we have
\begin{equation}\label{e:int-psi}
	\begin{aligned}
		\int_{- \pi N^{1/\rho} }^{ \pi N^{1/\rho} } \left|\psi(u/N^{1/\rho})\right|^n \D u &\le  2+  \int_{\{1\le |u|\le  \pi N^{1/\rho} \}} \left|\psi(u/N^{1/\rho})\right|^n \D u \\
		&\le 2+\int_\R \e^{-c|u|^{\rho-\ep}n/N} \D u\\
		&\le C (n/N)^{-\frac1{\rho-\ep}}. 
	\end{aligned}
\end{equation} 

Under the condition \eqref{e:con}, we may choose $\ep>0$ sufficiently small such that $2H-\frac1{\rho-\ep}>1$. To estimate $A_m$ given in \eqref{e:BNm}, we combine the estimate \eqref{e:int-psi} with the argument used in the proof of Proposition \ref{prop:exp-int} which leads to \eqref{e:estimation-moment},  and we can get for all $N\in \NN$, 
\begin{align*}
	A_m\le& \frac{(2m)!}{m!} C^m\int_{[0,1]_<^{2m}} \prod_{i=1}^m |t^*_{2i}-t^*_{2i-1}|^{2H-2-\frac1{\rho-\ep}}  \D \bbs\\
	\le &\frac{  (2m)! C^m}{m!\Gamma\left(m(2H-\frac1{\rho-\ep})+1\right)} 
\end{align*}
Therefore, we have uniformly in $N$, 
\[ \sum_{m=0}^\infty A_m \le \sum_{m=0}^\infty \frac{(2m)!C^m}{m!\Gamma\left(m(2H-\frac1{\rho-\ep})+1\right)} <\infty,\]
where the finiteness follows from Stirling's formula and $2H-\frac1{\rho-\ep}>1$. This implies
\begin{equation}\label{e:bd-partition}
	\sup_{N\in\NN} \E\left[\exp\left(N^{-\theta} \sum_{n=1}^N \om(n, S_n)\right)\right]<\infty. \end{equation}
Finally, the above analysis also yields
\begin{equation}\label{e:uniform-convergence}
	\lim_{M\to\infty} \sup_{N\in\NN} \sum_{m=M+1}^\infty \E\left[ \frac{1}{m!}\left|N^{-\theta} \sum_{n=1}^N \om(n, S_n)\right|^m\right]=0,
\end{equation}
which shall be used to prove the weak convergence of rescaled partition functions. 

Now we consider the It\^o-Skorohod case under the condition \eqref{e:con'}. Recall the partition function $\tilde Z_\om^{(N)}$ in \eqref{esptp'}.  Noting that $\mathbf S_m^{(N)}$ and $\mathbf S_n^{(N)}$ given in \eqref{e:S-m'} are orthogonal in $L^2(\Omega)$ if $m\neq n$, we have 
\[\E\left[\left|\sum_{m=0}^\infty\frac{1}{m!}{\mathbf S}_m^{(N)}\right|^2\right] =\sum_{m=0}^\infty \frac1{(m!)^2} \E\left[\left|\mathbf S_m^{(N)}\right|^2\right]. \]

Applying \eqref{e:E-wick-prod'}, one can calculate  $\frac1{(m!)^2}  \E\left[\left|\mathbf S_m^{(N)}\right|^2\right]$ to get the same upper bound uniformly in $N$ (up to a multiplicative  constant $C^m$) as $m!\|\hat \gg_m(\cdot; t,x)\|^2_{\mathscr H^{\otimes m}}$ (see \eqref{e:m-chaos} in Proposition \ref{prop:g-m} and see also Remark \ref{rem:bd-discrete}).
Thus, assuming \eqref{e:con'} we have 
\begin{equation}\label{e:uniform-convergence'}
	\lim_{M\to\infty} \sup_{N\in\NN} \E\left[\left|\sum_{m=M+1}^\infty\frac{1}{m!}{\mathbf S}_m^{(N)}\right|^2\right]=0.
\end{equation}

{\bf Acknowledgement}\quad  The authors would like to thank Rongfeng Sun and the referees for very helpful comments.

\appendix \section{Physical Wick product}\label{sec:wick}

\renewcommand{\theequation}{A.\arabic{equation}}

In order to expand the partition function  \eqref{partition2} in a proper way to obtain its weak convergence, we shall invoke the notion of \emph{physical} Wick product for general random variables (in contrast to the \emph{probabilistic} Wick product  defined via Wiener chaos in Section \ref{sec:malliavin}). The physical Wick product (also known as  Wick power, Wick polynomial or Wick renormalization) was introduced by Wick~\cite{wick1950} in  the study of quantum field theory.  We collect some facts on physical Wick products in this subsection. The reader is referred to \cite{at87,at06, giraitis1986multivariate, gjessing1993wick, lm16, surgailis1983} for more readable account. 

Let $\{X_i\}_{i\in \NN}$ be a family of real random variables with finite moments of all orders. The physical Wick product $:\!X_1\cdots X_n\!:$ is  defined recursively as follows. For $n=0$, we set $:~:=1$, and for $n\ge 1$,
\[\frac{\partial }{\partial X_i} :\! X_1\cdots X_n\!:=:\! X_1\cdots\hat X_i\cdots X_n\!:\,, \quad   \E\left[:\! X_1\cdots X_n\!:\right]=0,\]
where $\hat X_i$ means the absence of $X_i$ in the product.  For example, 
\begin{align*}
	:\! X_1\!:&=X_1-\E[X_1];\\
	:\!X_1 X_2\!:&=X_1X_2-X_1\E [X_2]-X_2\E [X_1]+2\E [X_1] \E [X_2]-\E[X_1X_2].
\end{align*}
We remark  that  different indices may refer to the same random variable. In this situation, as an example, we also write $:\! X^n Y^m\!: \deq :\! X_1\cdots X_nY_1\cdots Y_m\!:$,  if $X_1=\cdots=X_n=X$ and $Y_1=\cdots =Y_m=Y$. 

If we assume  $\E\left[e^{\beta \sum_{i=1}^n |X_i|}\right]<\infty$ for some $\beta>0$, the physical Wick product  can be equivalently defined by 
\begin{align}\label{wick}
	:\! X_1\cdots X_n\!:\deq \frac{\partial^n}{\partial z_1\cdots \partial z_n}  G(z_1, \cdots, z_n; X_1,\cdots, X_n)\bigg |_{ z=0}\, ,
\end{align}
where  \[G( z_1, \cdots, z_n;X_1, \cdots, X_n)\deq \frac{\e^{\sum_{i=1}^n z_iX_{i}}} {\E\left[\e^{\sum_{i=1}^n z_iX_{i}}\right]}\] is called the generating function (or Wick exponential).  If $\{X_i\}_{i\in \NN}$ is a centred Gaussian family, the generating function is simply $G( z, \bbX)=\e^{ z\cdot \bbX - z\cdot Q z/2}$ where $Q=(\E[X_iX_j])_{n\times n}$ is the covariance matrix of $\bbX=(X_1, \cdots, X_n)$, and the resulting Wick products are related to Hermite polynomials (see \eqref{e:hermit} in Section \ref{sec:malliavin}). In particular, for a Gaussian random vector $(X_1, \dots, X_n)$, the {\it physical} Wick product $:\!X_1\cdots X_n\!:$ coincides with {\it probabilistic} Wick product $X_1\diamond X_2 \diamond \cdots \diamond X_n$ defined in Section \ref{sec:malliavin}.

We collect some basic properties of physical Wick products. Clearly  $:\!X_1X_2\dots X_n\!:$ only involves the random variables $X_1, \dots, X_n$ and their joint moments up to order $n$, and it is \emph{symmetric} and \emph{multilinear} in $(X_1, \dots, X_n)$ (multilinear means  linear in terms of each $X_i, i=1, \dots, n$).  If two groups of random variables $\{X_1,X_2,\dots,X_m\}$ and $\{X_{m+1},\dots,X_n\}$ are  independent of  each other,  then $ :\!X_1\cdots X_n\!:=:\! X_1\cdots X_m\!::\! X_{m+1}\cdots X_n\!:$ (see \cite[eq. (2.4)]{at87}). We remind that  the physical Wick product is no longer associative, which is different from the ordinary product. For instance, assuming $\E[X]=0$, we have $:\!XXX\!: = X^3-3X\E[X^2] -\E[X^3]$ which is different from $:\!XY\!:|_{Y=:X^2:}=X^3-X\E[X^2]-\E[X^3]$.   Thus a physical Wick product  is a single term whose value is determined by the definition, and  cannot be viewed as a composition of two (or several) physical Wick products. For instance, all the physical Wick products $:\!X^3\!:, :\!X^2X\!:$ and $:\! X X^2\!:$ mean  the same $:\!XXX\!:$, and in particular  $:\! X X^2\!:\neq :\!XY\!:|_{Y=:X^2:}$.

For any finite index set $A=\{i_1, \dots, i_n\}\subset \NN$,  
we denote $:\!X_A\!:\deq :\!X_{i_1}\cdots X_{i_n}\!:$ and   similarly, we take the notation $X_A\deq\prod_{i\in A}X_i$ for the ordinary product.  We use $X^A$ to denote the set $\{X_i, i\in A\}$ of random variables. 

We recall some facts about cumulants. 
Let $\kappa(X^A)$ denote the joint cumulant of $X^A=\{X_i, i\in A\}$. Then 
\begin{equation}\label{e:cumulant}
	\kappa(X^A) =\sum_{V} (|V|-1)!(-1)^{|V|-1}\prod_{i=1}^{|V|} \E[X_{V_i}]
\end{equation}
and
\begin{equation}\label{e:expectation}
	\E[X_A] =\sum_{V} \prod_{i=1}^{|V|} \kappa(X^{V_i}),
\end{equation}
where the sum $\sum_V$ is taken over all partitions $V=\{V_1, \dots, V_k\}, k\ge 1$ of $A$, and $|V|=k$ is the partition number. For instance, 
\[\kappa(X_1) =\E[X_1], ~\kappa(X_1, X_2) = \E[X_1X_2]-\E[X_1]\E[X_2],\]
and 
\begin{align*}\kappa(X_1, X_2, X_3)=&\E[X_1X_2X_3]-\E[X_1X_2]\E[X_3]-\E[X_1X_3]\E[X_2]\\&-\E[X_2X_3]\E[X_1]+2\E[X_1]\E[X_2]\E[X_3].
\end{align*}

If $\{X_i\}_{i\in\NN}$ is \emph{Gaussian}, we have $\kappa(X^A)=0$ if $|A|\ge 3$. If we  assume further $\E[X_i]=0$, the formula \eqref{e:expectation} reduces to the Wick's theorem:
\begin{equation}\label{e:expectation'}
	\E[X_A] =\begin{cases}
		0,&  \text{ if  $|A|$ is odd},\vspace{0.2cm}\\\displaystyle \sum_{V} \prod_{\{i,j\}\in V} \E[X_iX_j],& \text{ if  $|A|$ is even},
	\end{cases}
\end{equation}
where the summation $\sum_{V}$ is taken over all \emph{pair partitions} $V=\{V_1, \dots, V_{|A|/2}\}$ of $A$.

Ordinary products and  physical Wick products are connected by the following formula (see \cite[Prop. 1]{surgailis1983} or \cite[Appendix B]{at06}),
\begin{equation}\label{wick}
	X_A=\sum_{B\subset A}:\!X_B\!:  \sum_{V}\prod_{i=1}^{|V|} \kappa(X^{V_i})=\sum_{B\subset A}:\!X_B\!: \E[X_{A\backslash B}],
\end{equation} 
and
\begin{equation}\label{wick'}
	:\!X_A\!:=\sum_{B\subset A}X_B  \sum_{V}(-1)^{|V|}\prod_{i=1}^{|V|}\kappa(X^{V_i}),
\end{equation}
where the sum $\sum_{B\subset A}$ is taken over all subsets $B\subset A$ including $B=\emptyset$,  and the sum $\sum_{V}$ is over all partitions $V=\{V_1, \dots, V_k\}, k\ge 1$ of the set $A\backslash B$. We use the convention $X_\emptyset=:\!X_\emptyset\!:=\kappa(X_{\emptyset})=1.$

The following formula (see \cite[Appendix B]{at06} or \cite[Lemma 4.5]{hs17}) will be used 
\begin{equation}\label{e:E-wick-prod}
	\E[:\! X_A\!: :\! X_B\!:] =\sum_V \prod_{i=1}^{|V|}\kappa (X^{V_i}),
\end{equation}
where the summation $\sum_V$ is taken over all partitions $V=\{V_1, \dots, V_k\}, k\ge1$ of $A\cup B$ satisfying $V_i\cap A\neq \emptyset\neq V_i\cap B$ for each $V_i$.  In particular, if we assume $\{X_i\}_{i\in\NN}$ is a centered \emph{Gaussian} family,  equations \eqref{e:E-wick-prod} and \eqref{e:expectation'} yield
\begin{equation}\label{e:E-wick-prod'}
	\E[:\! X_A\!: :\! X_B\!:] =\begin{cases}
		0, & \text{ if } |A|\neq |B|\vspace{0.2cm}\\ \displaystyle
		\sum_V \prod_{\{i,j\}\in V}\E[X_iX_j],  & \text{ if  } |A|=|B|,
	\end{cases}
\end{equation}
where the summation $\sum_{V}$ is taken over all \emph{pair partitions} $V=\{V_1, \dots, V_{|A|}\}$ of $A\cup B$ such that  $V_k=\{i,j\}$ with $i\in A, j\in B$ for $k=1, \dots, |A|$.

\section{Some preliminaries on convergence of probability measures}\label{sec:appB}
\renewcommand{\theequation}{B.\arabic{equation}}
The following result can be found in \cite[Theorem 3.2, Chapter 1]{Billingsley}. 
\begin{lemma}\label{lem:weak-con}
	Consider random vectors $Y_n^{(N)}, Y^{(N)}, Y_n$ and $Y$, such that $Y^{(N)}_n\wc Y_n$ as $N\to \infty$, $Y_n\wc Y$ as $n\to\infty$, and $Y_n^{(N)}$ converges in probability  to $Y^{(N)}$ uniformly in $N$ as $n\to\infty$, then we have $Y^{(N)}\wc Y$ as $N\to\infty$. That is, assuming
	\[\xymatrixcolsep{11pc} \xymatrix{
		Y_n^{(N)} \ar[d]_{d}^{N \to\infty} \ar[r]^{\mathrm{in\ probability,\ uniformly\ in\ } N }_{n \to\infty} & Y^{(N)} \\Y_n \ar[r]^{d}_{n\to\infty} & Y,
	}
	\]
	we have $Y^{(N)}\wc Y$ as $N\to\infty.$
\end{lemma}

\begin{lemma}\label{lem:weak-con''}
	Let $k\ge1$ be an integer.  If $(X_1^{(n)}, \dots, X_k^{(n)})\wc (X_1, \dots, X_k)$ as $n\to\infty$, we have \[f(X_1^{(n)}, \dots, X_k^{(n)}) \wc f(X_1, \dots, X_k)\]  for any continuous function $f$. 
\end{lemma}

\begin{lemma}\label{lem:weak-con'}
	Let $k\ge1$ be an integer.  Suppose $(X_1^{(n)}, \dots, X_k^{(n)})\wc (X_1, \dots, X_k)$ as $n\to\infty$ and assume that for any subset $A$ of $\{1, \dots, k\}$, $\left\{\prod_{i\in A}X_i^{(n)}\right\}_{n\in \NN}$ is uniformly integrable. Then $:\!X_1^{(n)}\cdots X_k^{(n)}\!:\wc :\!X_1\cdots X_k\!:$ as $n\to \infty$. 
\end{lemma}
\proof Lemma \ref{lem:weak-con''} yields the weak convergence of the ordinary product $\prod_{i\in A}X_i^{(n)}\wc  \prod_{i\in A}X_i$ for any $A\subset \{1, \dots, k\}$ as $n\to\infty$. 
By \eqref{wick'}, it suffices to prove the convergence of the cumulants \[\lim_{n\to\infty}\kappa(X_i^{(n)}, i\in A) =\kappa(X_i, i\in A).\]
This follows from the Skorohod representation theorem, equation \eqref{e:cumulant}, and the assumption of uniform integrability for the products of $X_1^{(n)},\dots, X_k^{(n)}$. 
\qed

\begin{theorem}\label{thm:cw}[Cram\'er-Wold Theorem] As $n\to \infty$,
	$(X_1^{(n)}, \dots, X_k^{(n)}) \wc (X_1, \dots, X_k)$ if and only if  
	\[\sum_{i=1}^k a_i X_i^{(n)} \wc \sum_{i=1}^k a_i X_i, ~ \text{ for all $(a_1, \dots, a_k)\in\R^k$.}
	\] 
\end{theorem}

\section{Miscellaneous results}
\label{sec:appC}

Let $f:[0,1]^m\to \R$ be an integrable function. For a fixed $N\in \NN$, let $t_i=i/N$ for $i=0,1,\dots, N$ and denote $I_i=(t_{i-1}, t_i]$.  Define 
\begin{equation}\label{e:fN}
	f_N(t_1, \dots, t_m)\deq   N^m \int_{I_{\ell_1}\times \cdots \times I_{\ell_m}} f(s_1, \dots, s_m) \D \bbs, 
\end{equation}
if $(t_1, \dots, t_m)\in I_{\ell_1}\times \dots \times I_{\ell_m}$  for some $\ell_i\in\llb N\rrb.$  That is, $f_N$ is the conditional expectation of $f$ with respect to the the $\sigma$-field generated by the rectangles of the form $I_{\ell_1}\times \dots \times I_{\ell_m} $. Then we have the following  inequality of Jensen type. 

\begin{lemma}\label{lem:jensen}
	Assume $H\in(1/2, 1]$ and suppose $f,g:[0,1]^{m}\to \R$  are integrable functions.  Then, for each $N\in \NN$,
	\begin{equation}\label{e:jensen}
		\begin{aligned}
			&\int_{[0,1]^{2m}}\prod_{i=1}^m |s_i-t_i|^{2H-2}|f_N(s_1,\dots,s_m)||g_N(t_1,\dots, t_m)|\D \bbs\D \bbt\\
			& \le C^m \int_{[0,1]^{2m}}\prod_{i=1}^m |s_i-t_i|^{2H-2}|f(s_1, \dots, s_m)||g(t_1, \dots, t_m)|\D \bbs\D \bbt,
		\end{aligned}
	\end{equation}
	for some constant $C$ depending on $H$ only. 
\end{lemma}
\proof	
We only prove the case $m=1$, and the other cases $m>1$ follows from an induction argument together with Fubini's theorem.  It suffices to prove \eqref{e:jensen} with $m=1$ for nonnegative functions, i.e., for $f, g\ge 0$, 
\begin{equation}\label{e:ineq-ij}
	\int_{I_i\times I_j} |s-t|^{2H-2}f_N(s)g_N(t)\D s\D t \le C \int_{I_i\times I_j} |s-t|^{2H-2}f(s)g(t) \D s\D t, ~ 1\le i, j \le N,
\end{equation} 
where recalling $I_i=(t_{i-1}, t_i]$.  We shall prove \eqref{e:ineq-ij} for simple nonnegative functions. The general case can be proved by a limiting argument and thus is omitted. 

For $i=j$, we assume that on $I_i$,
\[f(s)=\sum_{\ell=1}^k a_\ell \mathbf 1_{A_\ell}(s), ~ g(t)=\sum_{\ell=1}^k b_\ell \mathbf 1_{A_\ell}(t),\]
where $A_\ell = (t_{i-1}+(\ell-1)/Nk, t_{i-1}+\ell/Nk]$ with $\ell=1, \dots, k$ form a uniform partition of the interval $I_i$,  and $a_\ell, b_\ell$ are nonnegative numbers. Then, denoting $\bar a=\frac1k\sum_{\ell=1}^k a_\ell, \bar b=\frac1k\sum_{\ell=1}^k b_\ell$,
\[f_N(s)= \bar a, ~ g_N(t)=  \bar b, \text{ on } I_i. \]
The left-hand side of \eqref{e:ineq-ij} is 
\begin{equation}\label{e:lhs}
	\bar a\bar b \int_{I_i\times I_i} |s-t|^{2H-2} \D s \D t = \sum_{\ell=1}^k \sum_{m=1}^k a_\ell b_m \left( \frac1{k^2} \int_{I_i\times I_i} |s-t|^{2H-2} \D s \D t\right),
\end{equation}
and the right-hand side without the constant $C$ is
\begin{equation}\label{e:rhs}
	\sum_{\ell=1}^k \sum_{m=1}^k a_\ell b_m \left(\int_{A_\ell\times A_m}  |s-t|^{2H-2} \D s \D t \right). 
\end{equation}
Then, \eqref{e:ineq-ij} follows from \eqref{e:lhs}, \eqref{e:rhs} and the following inequality
\[  \frac1{k^2} \int_{I_i\times I_i} |s-t|^{2H-2} \D s \D t \le  C \min_{\ell,m\in\llb k \rrb}\int_{A_\ell\times A_m}  |s-t|^{2H-2} \D s \D t ,\]
where $C$ is a constant depending on $H$ only. Indeed, by change of variables, the above inequality is equivalent to 
\[\frac1{k^2} \int_0^{1}\int_0^{1} |s-t|^{2H-2} ds dt \le  C \min_{\ell,m\in\llb k \rrb}\int_{(m-1)/k}^{m/k}\int_{(\ell-1)/k}^{\ell/k}  |s-t|^{2H-2} \D s \D t,\]
which holds for $C=\int_0^1\int_0^1 |s-t|^{2H-2} \D s \D t$, noting that  for all $\ell, m\in \llb k \rrb$, we have $|s-t|\le 1$ on $[(\ell-1)/k, \ell/k] \times [(m-1)/k, m/k]$,  and hence $\int_{(m-1)/k}^{m/k}\int_{(\ell-1)/k}^{\ell/k} |s-t|^{2H-2} \D s \D t\ge 1/k^2$.
This proves \eqref{e:ineq-ij} for $i=j$. 

Now we prove \eqref{e:ineq-ij} for $i\neq j$. By symmetry, we only need to consider $j>i.$ Let
\[f(s)=\sum_{\ell=1}^k a_\ell \mathbf 1_{A_\ell}(s) \text{ on } I_i=(t_{i-1}, t_i],\]
and 
\[g(t)=\sum_{\ell=1}^k b_\ell \mathbf 1_{B_\ell}(t) \text{ on } I_j=(t_{j-1}, t_j],\]
where $a_\ell,b_\ell$ are nonnegative numbers, and $\{A_\ell, \ell=1, \dots, k\}$ and $\{B_\ell, \ell=1, \dots, k\}$ are uniform partitions of $I_i$ and $I_j$, respectively. Then,
\[f_N(s) = \bar a \text{ on } I_i \text{ and } g_N(t) = \bar b\text{ on } I_j,\]
where $\bar a=\frac{1}{k}\sum_{\ell=1}^k a_\ell$ and $\bar b=\frac{1}{k}\sum_{\ell=1}^k b_\ell$. The left-hand side of \eqref{e:ineq-ij} is 
\begin{equation}\label{e:lhs'}
	\bar a \bar b \int_{I_i\times I_j}|s-t|^{2H-2} \D s \D t = \sum_{\ell=1}^k\sum_{m=1}^k a_{\ell}b_{m}\left( \frac1{k^2} \int_{I_i\times I_j} |s-t|^{2H-2} \D s \D t\right),
\end{equation}
and the  right-hand side of \eqref{e:ineq-ij} without $C$ is
\begin{equation}\label{e:rhs'}
	\sum_{\ell=1}^k \sum_{m=1}^k a_\ell b_m \left(\int_{A_\ell\times B_m}  |s-t|^{2H-2} \D s \D t \right).
\end{equation}
To get \eqref{e:ineq-ij}, it suffices to show 
\begin{equation}\label{e:key-ine}
	\frac1{k^2} \int_{I_i\times I_j} |s-t|^{2H-2} \D s \D t \le C \min_{\ell, m\in\llb k\rrb} \int_{A_\ell\times B_m}  |s-t|^{2H-2} \D s \D t
\end{equation}
for some constant $C$ depending on $H$ only. Note that by change of variables, \eqref{e:key-ine} is equivalent to the following equality
\begin{equation}\label{e:key-ine'}
	\frac1{k^2} \int_{j-i}^{j-i+1} \int_{0}^{1} |s-t|^{2H-2} \D s \D t \le C \min_{\ell, m\in\llb k\rrb}   \int_{j-i+(m-1)/k}^{j-i+m/k} \int_{(\ell-1)/k}^{\ell/k} |s-t|^{2H-2} \D s \D t  
\end{equation}
We shall prove \eqref{e:key-ine'} for $j=i+1$ and $j\ge i+2$ separately. If $j=i+1$, \eqref{e:key-ine'} follows directly by choosing \[C= 2^{2-2H}\int_1^2  \int_0^1 |s-t|^{2H-2} \D s\D t,\] 
noting that
\[ \int_{1+(m-1)/k}^{1+m/k}  \int_{(\ell-1)/k}^{\ell/k} |s-t|^{2H-2} \D s \D t  \ge 2^{2H-2}/k^2, \text{ for all } \ell, m\in\llb k\rrb.\]

When $n\deq j-i\ge 2$, to prove \eqref{e:key-ine'}, it suffices to prove that the biggest integral among the integrals over the regions $\{(n+(m-1)/k, n+m/k]\times ((\ell-1)/k, \ell/k],\ell\in \llb k\rrb\}$  can be dominated by the smallest one, i.e., 
\[\int_{n}^{n+1/k}  \int_{1-1/k}^1 |s-t|^{2H-2} \D s \D t\le C \int_{n+1-1/k}^{n+1} \int_0^{1/k} |s-t|^{2H-2} \D s \D t\]
for some constant $C$ only depending on $H$. This is true because there exists a finite constant depending on $H$ only such that  
\[ |n-1|^{2H-2} \le 3^{2-2H} |n+1|^{2H-2} \text{ for all } n\ge 2.\]
This proves \eqref{e:key-ine'}, and hence completes the proof of \eqref{e:ineq-ij} for $i\neq j$.  \qed

%
%
%


The following Hardy-Littlewood-Sobolev inequality is taken from \cite[Lemma B.3 ]{balan2016intermittency} (see \cite{mmv01} for the one-dimensional version).
\begin{lemma}\label{lem:HLS}
	For $H\in(\frac12, 1]$, the following inequality holds:
	\begin{equation*}
		\int_{\R^m}\int_{\R^m}f(\bbt)f(\bbs)\prod_{i=1}^{m}|t_i-s_i|^{2H-2}\D \bbt\D\bbs\le  C_H^m \left( \int_{\R^m}\left|f(\bbt)\right|^{1/H}\D\bbt \right)^{2H},
	\end{equation*}
	where $C_H>0$ is a constant depending on $H$, and we denote $\bbt=(t_1,\dots,t_n)$ and $\bbs=(s_1,\dots,s_n)$.
\end{lemma}

The following result can be calculated by a direct calculation. 
\begin{lemma}\label{lem:ineq-gamma}
	Suppose $\alpha_i<1$ for  $i=1, \dots, m$ and let $\alpha=\sum_{i=1}^m\alpha_i$.
	Then 
	\[\int_{[0< r_1<\dots<r_m<r_{m+1}=t]}~ \prod_{i=1}^{m}( r_{i+1}-r_{i})^{-\alpha_i}\D \bbr = \frac{\prod_{i=1}^m\Gamma(1-\alpha_i)}{\Gamma(m-\alpha+1)}\, t^{m-\alpha},\]
	where $\Gamma(x)=\int_0^\infty t^{x-1} e^{-t}\D t$ is the Gamma function.
\end{lemma}




\end{document}